\documentclass[11pt, a4paper]{amsart}
\usepackage{amssymb,amsmath,amsfonts,latexsym,amsthm,mathtools}
\usepackage{enumitem}
\usepackage{comment}
\usepackage[margin=1.02in]{geometry}
\usepackage[colorlinks=true,linkcolor=blue,citecolor=blue]{hyperref}
\newtheorem{Theorem}{Theorem}[section]
\newtheorem{Lemma}[Theorem]{Lemma}
\newtheorem{Proposition}[Theorem]{Proposition}
\newtheorem{Corollary}[Theorem]{Corollary}
\newtheorem{theoremletter}{Theorem}

\theoremstyle{definition}
\newtheorem{Definition}[Theorem]{Definition}
\newtheorem{Example}[Theorem]{Example}
\newtheorem{Remark}[Theorem]{Remark}
\newtheorem{Question}[Theorem]{Question}
\numberwithin{equation}{section}

\def\Conf{\operatorname{Conf}}
\def\scal{\operatorname{scal}}
\def\Mo{\mathcal{M}}
\def\Po{{\mathfrak P}}
\def\pl{{\mathfrak p}_\ell}
\def\ql{{\mathfrak q}_\ell}
\def\m{{\mathfrak m}}
\def\hm{\hat{\mathfrak m}}
\def\R{{\mathbf R}}
\def\Q{{\mathbf Q}}
\def\C{{\mathbf C}}
\def\HH{{\mathbf H}}
\def\Hh{{\mathcal H}_\ell}
\def\kp{\varkappa}
\def\rank{\operatorname{rank}}
\def\Sim{\operatorname{Sim}}
\def\Isom{\operatorname{Isom}}
\def\Mob{\operatorname{M\ddot ob}}
\def\Symz{{\mathcal S}_\ell}
\def\im{\operatorname{im}}
\def\cent{\operatorname{cent}}
\providecommand{\MopsL}{\overline{\Mo}^{\,\mathrm{ps}}_\ell(S^n)}
\providecommand{\Xss}{X^{\mathrm{ss}}}
\providecommand{\Xst}{X^{\mathrm{s}}}
\providecommand{\Gmod}{\mathcal G_{n,\ell}(s)}
\providecommand{\Polys}{\mathcal P_\ell(s)}
\providecommand{\pt}{\operatorname{pt}}
\providecommand{\pAB}{{\mathfrak p}_{A|B}}
\let\phi=\varphi
\let\epsilon=\varepsilon
\dedicatory{Dedicated to the memory of Misha Kapovich}
\begin{document}

\title[Canonical metrics on punctured spheres, I]
      {Canonical metrics on punctured spheres, I:
        \\
the mass matrix and the polygon model}

\author{Boris Botvinnik}

\address{Department of Mathematics,
  University of Oregon, Eugene, OR 97403, USA}
\email{botvinn@uoregon.edu}


\begin{abstract}
A canonical metric $g_{\mathfrak{p}_{\ell}}$ on $S^n \setminus
\mathfrak{p}_{\ell}$, where
$\mathfrak{p}_{\ell}=(p_1,\ldots,p_{\ell})$ and $n \ge 3$, is a metric
conformal to the round metric that is complete, scalar-flat, and
asymptotically flat at each puncture $p_i$. For $n=3$ these metrics
are the Brill--Lindquist time-symmetric initial data for $\ell$ black
holes, and their ends carry ADM masses. In this paper we study the
masses, through the mass matrix, as functions on the conformal moduli
space $\mathcal{M}_{\ell}(S^n)$.

The mass matrix $\eta =(\eta_{ij})$ is a symmetric matrix with zero
diagonal and vanishing row sums (these matrices form the space
$\mathcal{H}_{\ell}$), and it determines the \emph{normalized ADM
mass} at the $i$-th end via $\hat{\mathfrak{m}}_i=\sum_{j\ne
  i}e^{\eta_{ij}}$. We prove that the mass map
$\Xi\colon\mathcal{M}_\ell(S^n)\to\mathcal{H}_{\ell}$ is a
homeomorphism onto its image.

On the other hand, the moduli space $\mathcal{M}_\ell(S^n)$ is
homeomorphic to the moduli space of those equilateral $\ell$-gons in
the hyperbolic space $\mathbb{H}^{n+1}$ whose sides, extended to
infinity, end at distinct points (the Gauss map). By a theorem of
Kapovich--Leeb--Millson and Balser every polystable configuration is
the Gauss image of such a polygon, and the Gauss map identifies the
space of all equilateral $\ell$-gons with a natural compactification
of $\mathcal{M}_\ell(S^n)$. In particular, the mass matrix has an
explicit description in terms of the pairwise vertex distances of the
polygon.

The mass matrix diverges as punctures collide.  The degeneration, its
renormalization, and the resulting compactification of the moduli
space are the subject of a sequel to this paper.
\end{abstract}
\vspace*{-10mm}

\maketitle
\vspace*{-12mm}

\renewcommand{\baselinestretch}{.7}
\setcounter{tocdepth}{1}
\begin{small}

\tableofcontents
\end{small}
\vspace*{-10mm}

\renewcommand{\baselinestretch}{1}
\section*{Introduction}\label{s0}

\subsection*{The geometric setting}
Let $h$ be the standard round metric of curvature one on $S^n$,
$n\ge3$, and let $\pl=(p_1,\dots,p_\ell)$ be an ordered configuration
of $\ell\ge3$ distinct points. We denote by
\[
\Conf_\ell(S^n) \coloneqq \{(p_1,\dots,p_\ell)\in(S^n)^\ell : p_i\ne p_j\text{ for }i\ne j\}
\]
the configuration space of $\ell$ labelled points, and by $\Mob(S^n)
\cong \mathrm O^\uparrow(1,n+1)$ the full M\"obius group acting on $S^n$. The
canonical metric on the punctured sphere
$S^n\setminus\{p_1,\dots,p_\ell\}$ is defined by
\[
g_{\pl} = \Phi_{\pl}^{4/(n-2)}h,\qquad
\Phi_{\pl}(x) = \sum_{i=1}^{\ell}c_i|x-p_i|^{-(n-2)},
\]
where the positive constants $c_i$ are chosen so that the metric $g_{\pl}$
depends only on the conformal equivalence class $[\pl]$ in the moduli
space
\[
\Mo_\ell(S^n) \coloneqq \Conf_\ell(S^n)/\Mob(S^n).
\]
Since $|x-p_i|^{-(n-2)}$ is, up to a constant factor, the Green's
function of the conformal Laplacian
$L_h=\frac{4(n-1)}{n-2}\Delta_h+\scal_h$ of the round sphere, the metric $g_{\pl}$ is complete,
scalar-flat, and asymptotically flat with $\ell$ ends. In dimension
$n=3$, these metrics represent the classical Brill--Lindquist
time-symmetric initial data for $\ell$ black holes; for $n>3$, they
provide the natural higher-dimensional generalizations.

The primary goal of this paper is to describe the ADM masses as
functions on the moduli space $\Mo_\ell(S^n)$. In the next paper
\cite{PartII}, we study their asymptotics and renormalization under
collisions of punctures.

\subsection*{Two normalizations in the light cone}
By identifying $S^n$ with the space of future-pointing null rays in
Minkowski space $\R^{1,n+1}$, the M\"obius group acts on $S^n$ via the
linear Lorentz group $\mathrm{O}^\uparrow(1,n+1) \cong
\Isom(\HH^{n+1})$. Under this identification, a configuration of
$\ell$ points corresponds to a collection of $\ell$ null
rays. Choosing a lift to specific null vectors is equivalent to
specifying a conformal factor at the points of $\pl$. This framework
naturally admits two distinct normalizations.
\begin{enumerate}[label=\textup{(\roman*)},leftmargin=*]
\item \emph{The mass normalization:} We set $c_i=e^{f_i}$, where $f$ is the unique solution to
\[
\eta_{ij} = -(n-2)\ln|p_i-p_j| + f_i + f_j\quad (i\ne j),
\qquad \sum_{j\ne i}\eta_{ij} = 0.
\]
The resulting matrix $\eta$ is a complete conformal invariant that encodes the normalized ADM
masses via $\hm_i = \frac{n-2}{4}\m_i = \sum_{j\ne i}e^{\eta_{ij}}$.
\item \emph{The polygon normalization:} We scale the null vectors so that the recursion
$x_{i+1} = e^{-s}x_i + w_i$ determines points on the unit hyperboloid in $\R^{1,n+1}$. The
points $x_i$ form the vertices of a closed equilateral $\ell$-gon in $\HH^{n+1}$ whose ideal edges
point toward the punctures.
\end{enumerate}

\subsection*{Main results}
Our first result is a normalized version of the rigidity theorem of
Crane and Short \cite{CS}, which states that two configurations are
M\"obius equivalent if and only if their absolute cross-ratios agree.
The mass matrix packages the cross-ratios in a normalization that ties
them to the ADM masses.

Denote by $\Hh$ the space of symmetric $\ell\times\ell$-matrices $(h_{ij})$
with zero diagonal and $\sum_{j\ne i}h_{ij} = 0$ for each
$i=1,\ldots,\ell$.
\begin{theoremletter}[After Crane--Short]\label{thm:mainA}
Let $\ell\ge3$ and $n\ge3$. The mass map $\Xi\colon
\Mo_\ell(S^n)\to\Hh$ is injective and proper, hence a homeomorphism
onto its image, which is closed in $\Hh$. A matrix $\eta\in\Hh$ lies in the image of
the mass map $\Xi$ if and only if the Gram matrix
\[
M(\eta)_{ij} = e^{-2\eta_{ij}/(n-2)}\quad (i\ne j), \qquad M(\eta)_{ii}=0
\]
has exactly one positive eigenvalue and at most $n+1$ negative eigenvalues.
\end{theoremletter}
The normalization condition $\sum_{j\ne i}\eta_{ij}=0$ is what makes
the mass map conformally invariant.

The moduli space $\Mo_\ell(S^n)$ is related to the moduli space of
equilateral $\ell$-gons in hyperbolic space $\HH^{n+1}$ modulo its
isometry group. Such an $\ell$-gon is a collection of points $x=(x_1,\ldots,x_{\ell})\in
(\HH^{n+1})^{\ell}$ such that the cyclic distances along the
connecting geodesics $\gamma_i$ satisfy $d(x_i,x_{i+1})=s$. The
infinite extension $\bar\gamma_i$ of each geodesic determines a unique
point on the sphere $S^n$, the boundary at infinity of
$\HH^{n+1}$. This is the Gauss map. It identifies $\Po_\ell(s)$ with the
moduli space
of ``polystable equally weighted configurations''
$\overline{\Mo}^{\,\mathrm{ps}}_\ell(S^n)$, a compactification of
the original moduli space $\Mo_\ell(S^n)$, see Section \ref{ss:bij}
for the definitions and Corollary~\ref{C:compact} for the details.
\begin{theoremletter}[After Kapovich--Leeb--Millson and Balser]\label{thm:mainB}
Let $s>0$. The hyperbolic Gauss map identifies the space $\Po_\ell(s)$
of closed $s$-equilateral $\ell$-gons in $\HH^{n+1}$ modulo isometry
with the equally weighted semistable quotient of \cite{KLM}, whose
points are the polystable orbits:
\[
\Gamma\colon \Po_\ell(s) \xrightarrow{\ \cong\ } \overline{\Mo}^{\,\mathrm{ps}}_\ell(S^n).
\]
Under this identification, on the locus of pairwise distinct
punctures, the ADM masses $\hm_i$ are explicit functions of the
pairwise hyperbolic vertex distances $d(x_a,x_b)$.
\end{theoremletter}

The existence of a polygon which maps to a prescribed configuration
$\pl=(p_1,\dots,p_\ell)$ via Gauss was proven in
Kapovich--Leeb--Millson \cite{KLM,KLM2} and Balser \cite{Ba}, and the
special case $n=2$ is in \cite{KMT}. We add the uniqueness of the
polygon in the stable case, from the fact that the composite step map
strictly decreases distances (Lemma~\ref{L:contract}), the
identification of the target with a compactification of
$\Mo_\ell(S^n)$ (Corollary~\ref{C:compact}), and the formula for the
masses (Proposition~\ref{P:massespoly}).

We also analyze the cases $\ell=4,5$ explicitly and show that for
$\ell\le n+2$ the image of the mass map $\Xi$ is star-shaped, so that
the moduli space $\Mo_\ell(S^n)$ is contractible, see also \cite{BS1}.

\subsection*{Plan of the paper}
Section~\ref{s:conv} sets notation and computes the ADM
masses. Section~\ref{s:massmatrix} constructs the mass matrix $\eta$
and establishes its canonicity and the tropical bounds relating it to
the masses. Section~\ref{s:lightcone} sets up the light-cone model and
proves Theorem~\ref{thm:mainA}. Section~\ref{s:polygons} constructs
the polygon model, proves Theorem~\ref{thm:mainB}, shows that
$\overline{\Mo}^{\,\mathrm{ps}}_\ell(S^n)$ is a compactification of
$\Mo_\ell(S^n)$, describes its boundary strata
(Proposition~\ref{P:psboundary}), and records the dictionary between initial data and
polygons. Section~\ref{s:small} treats $\ell=4,5$ and the shape-space
models.

\subsection*{Acknowledgement}
The author is deeply grateful to Kazuo Akutagawa, who initiated this
project in the early 2000s through his pioneering computations of the
ADM mass matrix in this setting.  This work is dedicated to the memory
of Misha Kapovich, who shared his mathematical insights with great
generosity. His influence took shape over many conversations
concerning the moduli spaces that are the subject of this paper.

\subsection*{AI disclosure}
Many of the constructions and first versions of the proofs in this paper
were developed in conversations with the AI model Claude Fable 5.1
(Anthropic), and the paper was revised with Claude Opus 5.5 (Anthropic).
The author substantially reworked all of this material into its present
form and takes responsibility for it.

\section{Conventions and the mass}\label{s:conv}

\subsection{Notation}
Let $h$ be the round metric of curvature $1$ on $S^n\subset\R^{n+1}$,
$n\ge3$, and $\phi(x,y)=|x-y|^2$ the squared chordal distance. Fix
$\pl=(p_1,\dots,p_\ell)\in\Conf_\ell(S^n)$, $\ell\ge3$, and write
\[
\kp:=n-2,\quad d_{ij}:=|p_i-p_j|,\quad \phi_{ij}:=\phi(p_i,p_j)=d_{ij}^2,\quad L_{ij}:=\ln d_{ij}\ (i\ne j),\quad L_{ii}:=0,
\]
\[
T_i:=\sum_{\alpha\ne i}L_{i\alpha},\qquad T:=\sum_{\alpha<\beta}L_{\alpha\beta}
=\tfrac12\sum_iT_i .
\]
For arbitrary positive constants $\tilde c_i$ consider
\begin{equation}\label{eq:phigen}
g=\tilde\Phi^{4/\kp}h,\qquad
\tilde\Phi(x)=\sum_i \tilde c_i\,\phi(x,p_i)^{-\kp/2}.
\end{equation}
The canonical metric is $g_{\pl}=\Phi_{\pl}^{4/\kp}h$ with
\begin{equation}\label{eq:phi}
\Phi_{\pl}(x)=\sum_i c_i\,\phi(x,p_i)^{-\kp/2},
\qquad
c_i=\Big\{\prod_{\alpha\ne i}\phi_{i\alpha}^{\frac{2}{\ell-1}}\!\!
\prod_{\substack{\alpha<\beta\\ \alpha,\beta\ne i}}\!\!
\phi_{\alpha\beta}^{-\frac{2}{(\ell-1)(\ell-2)}}\Big\}^{\kp/4}.
\end{equation}
In a chart adapted to the $i$-th end a metric \eqref{eq:phigen} takes the form
\begin{equation}\label{eq:bar_g}
\bar g_{ab}=\big(1+\m_i\rho^{2-n}+O(\rho^{1-n})\big)\delta_{ab},
\end{equation}
and $\m_i$ is the mass of that end; the ADM mass is $c(n)\m_i$ for a positive
dimensional constant $c(n)$ on which nothing below depends. We set
$\hm_i:=\frac{\kp}{4}\m_i$.

\begin{Proposition}[Akutagawa, unpublished]\label{P:mass}
For any positive constants $\tilde c_i$ the metric \eqref{eq:phigen}
is complete, scalar-flat and asymptotically flat with $\ell$ ends, and
with the convention \eqref{eq:bar_g},
\begin{equation}\label{eq:mass_formula}
\m_i=\frac{4}{\kp}\,\tilde c_i\sum_{j\ne i}\frac{\tilde c_j}{d_{ij}^{\,n-2}},
\qquad\text{equivalently}\qquad
\hm_i=\sum_{j\ne i}\xi_{ij},\quad \xi_{ij}:=\tilde c_i\tilde c_jd_{ij}^{-\kp}.
\end{equation}
\end{Proposition}

\begin{proof}
Take $i=1$, $p_1=(1,0)$ the north pole, and let
$\psi(y)=(1+|y|^2)^{-1}(|y|^2-1,2y)$ be the stereographic projection,
$p_j=\psi(y_j)$. Then
\[
\phi(\psi(y),\psi(z))=\frac{4|y-z|^2}{(1+|y|^2)(1+|z|^2)},\quad
\phi(\psi(y),p_1)=\frac{4}{1+|y|^2},\quad
\psi^*h=\frac{4\,\delta}{(1+|y|^2)^2},
\]
so $\tilde\Phi\circ\psi=2^{-\kp}(1+|y|^2)^{\kp/2}u$ with
$u(y)=\tilde c_1+\sum_{j\ge2}b_j|y-y_j|^{2-n}$ and
$b_j=\tilde c_j(1+|y_j|^2)^{\kp/2}=\tilde c_j(2/d_{1j})^{\kp}$, using
$1+|y_j|^2=4/d_{1j}^2$. Hence
\[
\psi^*g=(\tilde\Phi\circ\psi)^{4/\kp}\psi^*h=\tfrac14u^{4/\kp}\delta
=\tfrac14\tilde c_1^{4/\kp}\Big(1+\sum_{j\ge2}\tfrac{b_j}{\tilde c_1}|y-y_j|^{2-n}\Big)^{4/\kp}\delta .
\]
The function $u$ is positive and harmonic outside of the punctures
$y_j$, so $g$ is scalar-flat. Near each $y_j$ the metric is
asymptotically flat in an inverted coordinate, and near $y=\infty$ in a
rescaling of $y$ itself, as the following computation shows for that end. We set $y=\alpha
z$ with $\tfrac14\tilde c_1^{4/\kp}\alpha^2=1$, i.e.\ $\alpha=2\tilde
c_1^{-2/\kp}$, and $\rho=|z|$, and obtain:
\[
\bar g_{ab}=\big[1+\frac4\kp\frac{\alpha^{2-n}}{\tilde c_1}(\sum_jb_j)\rho^{2-n}
  +O(\rho^{1-n})\big]\delta_{ab},
\]
and $\alpha^{-\kp}=2^{-\kp}\tilde c_1^2$, whence
$\m_1=\frac4\kp \tilde c_1\sum_{j\ge2}\tilde c_jd_{1j}^{-\kp}$.
\end{proof}
We do not yet impose a normalization on these constants; the canonical
constants \eqref{eq:phi} are introduced in the next section.

\section{The mass matrix}\label{s:massmatrix}

\subsection{Canonical mass matrix \texorpdfstring{$\eta$}{eta}}
We rely on the following simple linear algebra lemma:
\begin{Lemma}\label{L:centering}
Let $\ell\ge3$ and let $(L_{ij})$ be a symmetric $\ell\times\ell$-matrix
with zero diagonal. The system
\begin{equation}\label{eq:eta_sys}
\eta_{ij}=-\kp L_{ij}+f_i+f_j\ (i\ne j),\qquad \sum_{j\ne i}\eta_{ij}=0\ (1\le i\le\ell)
\end{equation}
has the unique solution
\begin{equation}\label{eq:f_sol}
f_i=\frac{\kp}{\ell-2}\Big(T_i-\frac{T}{\ell-1}\Big).
\end{equation}
For $L_{ij}=\ln|p_i-p_j|$ one has $e^{f_i}=c_i$, the constant of \eqref{eq:phi}.
\end{Lemma}

\begin{proof}
With $F=\sum_\alpha f_\alpha$, \eqref{eq:eta_sys} reads
$(\ell-2)f_i+F=\kp T_i$; summing gives $(2\ell-2)F=2\kp T$,
i.e.\ $F=\kp T/(\ell-1)$, whence \eqref{eq:f_sol}, the system being
nonsingular for $\ell\ge3$. For the last claim take logarithms in
\eqref{eq:phi} and use $\ln\phi_{ij}=2L_{ij}$:
\[
\ln c_i=\frac{\kp}{4}\Big[\frac{4T_i}{\ell-1}-\frac{4(T-T_i)}{(\ell-1)(\ell-2)}\Big]
=\frac{\kp[(\ell-1)T_i-T]}{(\ell-1)(\ell-2)}
=\frac{\kp}{\ell-2}\Big(T_i-\frac{T}{\ell-1}\Big).\qedhere
\]
\end{proof}
We write $\cent(L)\in\R^\ell$ for the solution \eqref{eq:f_sol} of
\eqref{eq:eta_sys}, the \emph{centering vector} of $L$; it depends linearly
on $L$.

\begin{Definition}\label{D:H}
  Let $\Hh$ be the following space of $\ell\times\ell$-matrices
  \[
\Hh:=\{\eta\in\R^{\ell\times\ell}:\eta_{ij}=\eta_{ji},\ \eta_{ii}=0,\
\textstyle\sum_j\eta_{ij}=0\}\cong\R^{\ell(\ell-3)/2},
\]
with the inner product
$\langle\eta,\eta'\rangle=\sum_{i<j}\eta_{ij}\eta'_{ij}$. For $\pl\in
\Conf_\ell(S^n)$ let the \emph{mass matrix} $\eta(\pl)\in\Hh$ be given
by Lemma~\ref{L:centering}; thus $\eta_{ij}=\ln\xi_{ij}$ as in
\eqref{eq:mass_formula} with the canonical constants. Set
$\nu:=\eta/\kp$.
\end{Definition}

\begin{Proposition}\label{P:invariance}
Let $\gamma\in\Mob(S^n)$, $\pl\in\Conf_\ell(S^n)$, then
$\eta(\gamma\pl)=\eta(\pl)$, so the mass matrix $\eta$ descends to the
\emph{mass map} $\Xi:\Mo_\ell(S^n)\to\Hh$. Moreover $\nu$ is
independent of $n$, and $\hm_i=\sum_{j\ne i}e^{\eta_{ij}}=\sum_{j\ne
  i}e^{(n-2)\nu_{ij}}$.
\end{Proposition}

\begin{proof}
We have $\phi(\gamma p_i,\gamma p_j)=|\gamma'_{p_i}||\gamma'_{p_j}|\phi_{ij}$, so
$L^\gamma_{ij}=L_{ij}+\tfrac12(\lambda_i+\lambda_j)$ with
$\lambda_i=\ln|\gamma'_{p_i}|$: the conformal group acts on the matrix of
log-distances by addition of rank one symmetric terms, which is exactly the ambiguity
divided out in \eqref{eq:eta_sys}. The vector $(f_i+\tfrac\kp2\lambda_i)$ solves the system for
$L^\gamma$ and yields the same $\eta$; by uniqueness it is the solution. As a
by-product $c_i(\gamma\pl)=|\gamma'_{p_i}|^{\kp/2}c_i(\pl)$, which is the covariance
needed for $g_{\pl}$ to descend. The rest is immediate from \eqref{eq:f_sol} and \eqref{eq:mass_formula}.
\end{proof}

Lemma~\ref{L:centering} can be reformulated in terms of an
orthogonal projection; we use this reformulation in the proofs of
Theorem~\ref{T:image} and Proposition~\ref{P:massespoly}.

\begin{Lemma}\label{L:proj}
Write $\Symz$ for the space of symmetric $\ell\times\ell$ matrices
with zero diagonal, carrying the inner product of
Definition~\ref{D:H}, and put
\[
\varpi:\R^\ell\longrightarrow\Symz,\qquad \varpi(a)_{ij}=a_i+a_j\ \ (i\ne j),
\qquad \Lambda_\ell:=\im\varpi .
\]
For $\ell\ge3$ the map $\varpi$ is injective, so
$\dim\Lambda_\ell=\ell$; moreover
\[
\Symz=\Hh\oplus^\perp\Lambda_\ell ,
\]
and the map $\Pi:L\mapsto\eta$ of \eqref{eq:eta_sys}--\eqref{eq:f_sol}
equals $-\kp P$, where $P$ is the orthogonal projection of $\Symz$
onto $\Hh$. In particular $\Pi$ is linear with kernel $\Lambda_\ell$,
$\|\Pi(S)\|\le\kp\|S\|$, and $\cent(\varpi(a))=\kp a$.
\end{Lemma}
\begin{proof}
For $\eta\in\Hh$ and $a\in\R^\ell$,
\[
\langle\eta,\varpi(a)\rangle=\tfrac12\sum_{i\ne j}\eta_{ij}(a_i+a_j)
=\sum_i a_i\sum_{j\ne i}\eta_{ij}=0 ,
\]
so $\Hh\perp\Lambda_\ell$; and
$\dim\Hh+\dim\Lambda_\ell=\tfrac{\ell(\ell-3)}2+\ell=\binom\ell2=\dim\Symz$,
so the sum is the whole space. Next $\Pi(\varpi(a))=0$: trying
$f_i=\kp a_i+c$ in \eqref{eq:eta_sys} gives $\eta_{ij}=2c$, whose row
sums are $2c(\ell-1)$, so $c=0$ and, by the uniqueness in
Lemma~\ref{L:centering}, $\eta=0$ and $\cent(\varpi(a))=\kp
a$. Finally $\Pi|_{\Hh}=-\kp\,\mathrm{id}$, since for $L\in\Hh$ the
choice $f=0$ already solves \eqref{eq:eta_sys}. A linear map equal to
$-\kp\,\mathrm{id}$ on $\Hh$ and to $0$ on $\Lambda_\ell=\Hh^\perp$ is
$-\kp P$.
\end{proof}
\begin{Remark}
  \label{R:lambda}
The space $\Lambda_\ell$ warrants some explanation. As shown in the
proof of Proposition~\ref{P:invariance}, the conformal group
translates the matrix of log-distances according to the
formula:
\begin{equation}\label{eq:gamma}
  L^\gamma=L+\tfrac12\varpi(\lambda),\qquad\lambda=\lambda(\gamma,\pl)\in\R^\ell,\quad
  \lambda_i=\ln|\gamma'_{p_i}| .
\end{equation}
In the light-cone language of Remark~\ref{R:lift}(1), this is the
cocycle $\lambda_i=-\ln\sigma_\gamma(p_i)$, which measures the
deviation of the section $u_0=1$ from equivariance at the $\ell$
punctures. We emphasize that the vector $\lambda$ depends on both the
configuration $\pl$ and the transformation $\gamma$: the formula
describes an action on pairs $(L,\pl)$, rather than just a linear
action on $L$. Thus $\Lambda_\ell$ records the ambiguity introduced by
the conformal group, although the group does not act linearly on $L$.
\end{Remark}

\subsection{Canonicity of the metric and the mass map}

\begin{Corollary}\label{T:unique}
Let $\pl\in\Conf_\ell(S^n)$ and let $\tilde c_1,\dots,\tilde c_\ell$
be any positive constants such that $\tilde\xi_{ij}=\tilde c_i\tilde
c_jd_{ij}^{-\kp}$ satisfies $\prod_{j\ne i}\tilde\xi_{ij}=1$ for all
$i$. Then $\tilde c_i=c_i$, $i=1,\ldots,\ell$.
\end{Corollary}

\begin{proof}
Put $\lambda_i=\tilde c_i/c_i$, so that
$\tilde\xi_{ij}=\lambda_i\lambda_j\xi_{ij}$ with $\prod_{j\ne
i}\xi_{ij}=1$. With $\mu_i=\ln\lambda_i$, $M=\sum\mu_\alpha$,
\[
0=\sum_{j\ne i}\ln\tilde\xi_{ij}=(\ell-1)\mu_i+(M-\mu_i)=(\ell-2)\mu_i+M,
\]
so all $\mu_i$ are equal to some $\mu$ and then $(2\ell-2)\mu=0$.
\end{proof}
Thus the row-product normalization alone pins down the constants, and
the resulting $c_i$ are automatically conformally covariant
(Proposition~\ref{P:invariance}). Covariance alone would not suffice:
for $\ell\ge4$ the group $\Mob(S^n)$ has invariants, so covariant
constants $\tilde c_i(\pl)$ are determined only up to multiplication
by arbitrary positive functions on $\Mo_\ell(S^n)$, an
infinite-dimensional ambiguity.

\subsection{Lower bounds for the masses}

\begin{Lemma}\label{L:tropical}
For $\eta\in\Hh$ set $\|\eta_{i\cdot}\|_1=\sum_{j\ne i}|\eta_{ij}|$,
$\eta_i^{\max}=\max_{j\ne i}\eta_{ij}$. Then
\begin{equation}\label{eq:tropical}
\eta_i^{\max}\le\ln\hm_i\le\eta_i^{\max}+\ln(\ell-1), \qquad
\frac{\|\eta_{i\cdot}\|_1}{2(\ell-1)}\le\eta_i^{\max}\le\frac{\|\eta_{i\cdot}\|_1}{2}.
\end{equation}
\end{Lemma}

\begin{proof}
The first pair is
$e^{\eta_i^{\max}}\le\sum_je^{\eta_{ij}}\le(\ell-1)e^{\eta_i^{\max}}$.
For the second, $\sum_j\eta_{ij}=0$ gives
$\sum_j\eta_{ij}^+=\sum_j\eta_{ij}^-=\tfrac12\|\eta_{i\cdot}\|_1$,
whence $\eta_i^{\max}\ge\frac{1}{\ell-1}\sum_j\eta^+_{ij}$ and
$\eta_i^{\max}\le\sum_j\eta^+_{ij}$. Here $\eta_{ij}^{\pm}$ denote the
positive and negative parts of $\eta_{ij}$.
\end{proof}

\begin{Remark}
Lemma~\ref{L:tropical} implies that $\max_i\ln\hm_i\to\infty$ if and
only if $\|\eta\|\to\infty$. The normalization $\sum_j\eta_{ij}=0$ is
essential here: without it a row of $\eta$ could have all entries large
and negative, so that $\|\eta_{i\cdot}\|_1\to\infty$ while $\hm_i$ stays
bounded.
\end{Remark}

\begin{Corollary}\label{C:bound}
$\hm_i\ge\ell-1$, with equality iff the $i$-th row of $\eta$ vanishes;
  equality for all $i$ iff $\eta=0$ iff $d_{ij}=t_it_j$ for some
  $t_i>0$, a conformally invariant condition.
\end{Corollary}
\begin{proof}
By the arithmetic--geometric mean inequality and $\prod_{j\ne
  i}\xi_{ij}=1$,
\[
\hm_i\ge(\ell-1)\Big(\prod_{j\ne i}\xi_{ij}\Big)^{1/(\ell-1)}=\ell-1.
\]
Equality holds if and only if $\xi_{ij}\equiv1$ along the row. If this
holds for every $i$ then $\eta=0$, i.e.\ $\kp L_{ij}=f_i+f_j$,
i.e.\ $d_{ij}=t_it_j$ with $t_i=e^{f_i/\kp}$; and $L$ changes by
$\tfrac12(\lambda_i+\lambda_j)$ under $\Mob(S^n)$, so $t_i\mapsto
e^{\lambda_i/2}t_i$.
\end{proof}
\begin{Remark}\label{R:balancing}
In a stereographic chart, $\Phi=\sum_ic_i|x-p_i|^{-\kp}$ is the
potential of point charges $c_i$, and $\xi_{ij}=c_ic_j\,d_{ij}^{-\kp}$
is the pair energy of the charges at $p_i$ and $p_j$; the mass
$\hm_i=\sum_{j\ne i}\xi_{ij}$ is the energy of the $i$-th charge in
the field of the others. Changing the charges,
$c_i\mapsto\lambda_ic_i$, rescales the symmetric matrix $(\xi_{ij})$
by a diagonal congruence, $\xi\mapsto D\xi D$, and the normalization
\eqref{eq:eta_sys} selects the unique such rescaling with all row
products equal to one (Corollary~\ref{T:unique}). This is the
logarithmic, and therefore linear, version of the matrix balancing
problem: Sinkhorn's theorem \cite{Si64} rescales a positive matrix by
$D\xi D$ so that all row \emph{sums} equal one, which is nonlinear and
requires iteration, whereas row \emph{products} are fixed by the
linear system \eqref{eq:eta_sys}. The additive ambiguity
$\ln\xi_{ij}=\eta_{ij}+\mu_i+\mu_j$ is the choice of conformal gauge, and
$\eta$ is what is left when it is fixed canonically.
\end{Remark}

\section{The light-cone model}\label{s:lightcone}
\subsection{Homogeneous coordinates}
Let $\R^{1,N}$ be $\R^{N+1}$ with the Lorentz form
\[
\langle u,v\rangle=-u_0v_0+\sum_{a\ge1}u_av_a,
\]
and let $ \mathcal N^+=\{u\ne0:\langle u,u\rangle=0,\ u_0>0\} $ be the
future null cone in $\R^{1,n+1}$. Every null line meets $\mathcal N^+$
in one ray, so we can identify the sphere $S^n$ with the quotient
$\mathcal N^+/\R_{>0}$, which is the projective quadric. Then for a
point $p\in S^n$ the vector $\hat p:=(1,p)$ is a choice of
\emph{homogeneous coordinates}, namely the section $u_0=1$. We have an
obvious formula:
\begin{equation}\label{eq:chordal}
\langle\hat p,\hat q\rangle=-1+\langle p,q\rangle=-\tfrac12|p-q|^2 .
\end{equation}

\begin{Remark}
  \label{R:lift}
We have two consequences which we use throughout this paper.
\begin{enumerate}
\item[(1)] {\bf The group $\Mob(S^n)$ acts on rays, not on the lifts.}
  If $\Gamma\in O^\uparrow(1,n+1)$ induces the transformation
  $\gamma\in\Mob(S^n)$, then $\Gamma(\hat
  p)=\sigma_\gamma(p)\widehat{\gamma(p)}$ with $\sigma_\gamma(p)>0$,
  and $\sigma$ is a cocycle, (i.e.,
  $\sigma_{\gamma_1\gamma_2}(p)=\sigma_{\gamma_1}(\gamma_2p)\sigma_{\gamma_2}(p)$).
  Comparing with \eqref{eq:chordal} gives $\phi(\gamma p,\gamma
  q)=|\gamma'_p||\gamma'_q|\phi(p,q)$ and
  $\sigma_\gamma(p)=|\gamma'_p|^{-1}$: this shows that the section $u_0=1$ is
  not $\Mob(S^n)$-equivariant.
\item[(2)] {\bf Sections are conformal metrics.}  For
  $u(p)=\omega(p)\hat p$, $\omega>0$,
\[
\langle du,du\rangle=d\omega^2\langle\hat p,\hat p\rangle
+2\omega d\omega\langle\hat p,(0,dp)\rangle+\omega^2|dp|^2=\omega^2h,
\]
since $\hat p$ is null and $\langle\hat p,(0,dp)\rangle=p\cdot
dp=0$. So sections of $\mathcal N^+\to S^n$ are exactly the metrics in
the conformal class; the section cut by $\langle u,x\rangle=-1$ has
curvature $+1$, $0$ or $-1$ according as $x$ is timelike, null or
spacelike, the null case being stereographic projection. Choosing
scalings $t_i>0$ at $\ell$ punctures is choosing the values of a
conformal factor there.
\end{enumerate}
\end{Remark}

\subsection{The Gram matrix of the mass normalization}

\begin{Definition}\label{D:M}
Let $\eta\in\Hh$. The Gram matrix $M(\eta)$ of $\eta$ is defined as
$M(\eta)_{ij}=e^{-2\eta_{ij}/\kp}$ ($i\ne j$), $M(\eta)_{ii}=0$.
\end{Definition}

\begin{Lemma}\label{L:gram}
Let $\pl$ have mass matrix $\eta$ and centering $f$, and put
$v_i:=e^{-f_i/\kp}\hat p_i$. Then the $v_i$ are future null and
$\langle v_i,v_j\rangle=-\tfrac12M(\eta)_{ij}$.
\end{Lemma}

\begin{proof}
By \eqref{eq:eta_sys}, $L_{ij}=(f_i+f_j-\eta_{ij})/\kp$, so
$d_{ij}^2=e^{2(f_i+f_j)/\kp}M(\eta)_{ij}$, and by \eqref{eq:chordal},
$\langle v_i,v_j\rangle=e^{-(f_i+f_j)/\kp}(-\tfrac12d_{ij}^2)=-\tfrac12M(\eta)_{ij}$.
\end{proof}

Before the theorem we need two facts that make the mass map $\Xi$ a
homeomorphism and not merely an injection. Both are statements about
the Gram matrix $M(\eta)$, and they hold for all $\ell$ and $n$.

\begin{Lemma}\label{L:crossratio}
For distinct indices $i,j,k,l$,
\[
\frac{M(\eta)_{ij}M(\eta)_{kl}}{M(\eta)_{ik}M(\eta)_{jl}}
=\exp\Big(-\tfrac2\kp\big(\eta_{ij}+\eta_{kl}-\eta_{ik}-\eta_{jl}\big)\Big)
=\Big(\frac{d_{ij}\,d_{kl}}{d_{ik}\,d_{jl}}\Big)^{2} .
\]
Thus a cross-ratio of the configuration is a  continuous
function depending only on the mass matrix $\eta$.
\end{Lemma}

\begin{proof}
The first equality is Definition~\ref{D:M}. For the second,
\eqref{eq:eta_sys} gives $\eta_{ij}=-\kp L_{ij}+f_i+f_j$, and the
$f$'s cancel in the alternating combination, so
$\eta_{ij}+\eta_{kl}-\eta_{ik}-\eta_{jl}=-\kp\ln\frac{d_{ij}d_{kl}}{d_{ik}d_{jl}}$.
Taking exponentials completes the proof.
\end{proof}

\begin{Proposition}\label{P:proper}
The mass map $\Xi$ is proper: $\Xi^{-1}(K)$ is compact for every
compact $K\subseteq\Hh$. Consequently $\Xi$ is a closed map, and its
image is closed in $\Hh$.
\end{Proposition}

\begin{proof}
For $\ell=3$ the space $\Mo_3(S^n)$ is a point, so assume
$\ell\ge4$. Put $C=\max_{K}\max_{i\ne j}|\eta_{ij}|$ and
$A=e^{4C/\kp}\ge1$. If $\Xi([\pl])\in K$ then
$|\eta_{ij}+\eta_{kl}-\eta_{ik}-\eta_{jl}|\le4C$, so by
Lemma~\ref{L:crossratio} every absolute cross-ratio of $\pl$ lies in
$[A^{-1},A]$.

The group $\Mob(S^n)$ acts transitively on ordered triples of distinct
points, so we may normalize $p_1,p_2,p_3$ to be the vertices of an
equilateral triangle on a great circle,
$d_{12}=d_{13}=d_{23}=\sqrt3$. For $i\ge4$ the cross-ratios of
$(1,i,2,3)$ and of $(1,i,3,2)$ are $d_{1i}/d_{3i}$ and
$d_{1i}/d_{2i}$, so both ratios lie in $[A^{-1},A]$. Since
$d_{1i}+d_{2i}\ge d_{12}=\sqrt3$ by the triangle inequality for the
chordal distance, $\max(d_{1i},d_{2i})\ge\sqrt3/2$, hence
$d_{1i}\ge\sqrt3/(2A)$, and then $d_{2i},d_{3i}\ge d_{1i}/A$. So
$d_{1i},\ d_{2i},\ d_{3i}\ \ge\ \frac{\sqrt3}{2A^2}$ for $i\ge4$.
For $4\le i<j$ the cross-ratio of $(i,j,1,2)$ is
$\sqrt3\,d_{ij}/(d_{1i}d_{2j})$, and we obtain
\[
\begin{array}{c}
d_{ij}\ \ge\ \frac{1}{A\sqrt3}\Big(\frac{\sqrt3}{2A^2}\Big)^{2}=\frac{\sqrt3}{4A^{5}} .
\end{array}
\]
Thus all pairwise distances are bounded below by
$\varrho:=\sqrt3/(4A^{5})>0$, a bound depending only on $K$.  Now let
$[\pl^{(m)}]$ be a sequence in $\Xi^{-1}(K)$ and choose
representatives normalized as above. They lie in the compact subset
$\{(p_1,\dots,p_\ell)\in(S^n)^\ell:d_{ij}\ge\varrho\}\subset
(S^n)^\ell$.

Hence a subsequence of the representatives converges to a
configuration $\pl^\infty$ whose points are still distinct. Since
$\eta$ is continuous on $\Conf_\ell(S^n)$ and $K$ is closed,
$\eta(\pl^\infty)=\lim\eta(\pl^{(m)})\in K$, so
$[\pl^\infty]\in\Xi^{-1}(K)$ and the subsequence converges to it in
$\Mo_\ell(S^n)$. The action of $\Mob(S^n)$ on $\Conf_\ell(S^n)$,
$\ell\ge3$, is proper: it is transitive on $\Conf_3(S^n)$ with compact
stabilizer $\mathrm{O}(n-1)$, so the action on $\Conf_3(S^n)$ is
proper, and properness pulls back along the equivariant projection
$\Conf_\ell(S^n)\to\Conf_3(S^n)$. It follows that $\Mo_\ell(S^n)$ is
Hausdorff, locally compact and metrizable, and we obtain that
$\Xi^{-1}(K)$ is compact.

A proper continuous map into a locally compact Hausdorff space is
closed, so in particular its image is closed.
\end{proof}

\begin{Theorem}[Theorem~\ref{thm:mainA}, cf. \cite{CS}]\label{T:image}
Let $\ell\ge3$, $n\ge3$.
\begin{enumerate}
\item[{\rm (1)}] The mass map $\Xi:\Mo_\ell(S^n)\to\Hh$ is a homeomorphism onto its
  image, which is closed; in particular the mass matrix is a complete conformal
  invariant of a labelled configuration.
\item[{\rm (2)}] A matrix $\eta\in\Hh$ lies in the image iff $M(\eta)$
  has exactly one positive and at most $n+1$ negative eigenvalues.
\end{enumerate}
\end{Theorem}

\begin{proof}
(1) If $\eta(\pl)=\eta(\ql)$, the associated null vectors $v_i,w_i$ of
  Lemma~\ref{L:gram} have the same Gram matrix. Since the
  configurations contain at least two distinct points, the span of the
  $v_i$ contains a timelike vector and is therefore nondegenerate; the
  same holds for the span of the $w_i$. Thus the correspondence
  $v_i\mapsto w_i$ is a well-defined isometry of these spans and
  extends, by Witt's theorem, to a transformation $\Gamma\in
  O(1,n+1)$. As $v_1$ and $w_1=\Gamma v_1$ are both future-pointing,
  $\Gamma$ preserves the future cone, and the induced
  $\gamma\in\Mob(S^n)$ carries $p_i$ to $q_i$. So $\Xi$ is injective.
  By Proposition~\ref{P:proper} it is also continuous and closed, and a
  closed continuous injection is a homeomorphism onto its image, which
  is closed.

(2) \emph{Necessity.} The matrix $-\tfrac12M(\eta)$ is a Gram matrix
  of vectors in $\R^{1,n+1}$, whose inertia is $(n+1,1)$; hence it has
  at most $n+1$ positive and at most one negative eigenvalue, so $M$
  has at most one positive and at most $n+1$ negative. As
  $\operatorname{tr}M=0$ and $M\ne0$, it has at least one of each,
  hence exactly one positive.

\emph{Sufficiency.} If $M$ has one positive and $q\le n+1$ negative
eigenvalues then $-\tfrac12M$ has inertia $(q,1)$ and is realized as a
Gram matrix of nonzero vectors $v_i\in\R^{1,n+1}$; $M_{ii}=0$ makes them null,
while $M_{ij}>0$ gives $\langle v_i,v_j\rangle<0$ for $i\ne j$. Hence no two
$v_i$ are proportional, and all lie in the same component of the null cone; after an
overall time reversal if necessary, we may take that component to be the future one.
Writing $v_i=t_i\hat p_i$ with $t_i>0$ and $p_i$ distinct, \eqref{eq:chordal} gives
$d_{ij}^2=t_i^{-1}t_j^{-1}M_{ij}$, i.e.\ $$-\kp
L_{ij}=\eta_{ij}+\tfrac{\kp}{2}\varpi(\ln t)_{ij},$$ so $L\in
-\eta/\kp+\Lambda_\ell$, and Lemma~\ref{L:proj} gives
$\eta((p_i))=\Pi(L)=-\kp P(L)=\eta$.
\end{proof}
\begin{Remark}\label{R:CS}
Geometrically, Theorem~\ref{T:image}(1) recovers the
point-configuration rigidity theorem of Crane and Short \cite{CS},
which answered a question posed by Beardon and Minda \cite{BM}: two
collections of distinct points in $\R^N_\infty$ are M\"obius equivalent
if and only if their absolute cross-ratios of ordered $4$-tuples
agree. We adapt their fundamental argument in our proof above. In this
framework, points on the sphere correspond to future-pointing null
lines in the hyperboloid model of $\HH^{N+1}$. This identification
allows the cross-ratio to be expressed via Lorentz products $\langle
v_i,v_j\rangle$. The rigidity conclusion then relies on the algebraic
fact that two families of vectors with identical Lorentz Gram matrices
spanning a time-like or space-like subspace differ only by a Lorentz
transformation \cite[Prop.~7.1]{CS}. Our degeneracy condition mirrors
theirs. Specifically, the Lorentz transformation fails to be unique
precisely when the vectors $v_i$ span a proper subspace of
$\R^{1,n+1}$; geometrically, this occurs if and only if the punctures
lie on a proper round subsphere. We exploit this structure by
introducing a canonical normalization. By repackaging the cross-ratio
invariant as a mass matrix $\eta\in\Hh$ with a zero diagonal and zero
row sums, we explicitly tie the configuration to its ADM masses. This
allows us to formulate Theorem~\ref{T:image}(2) directly as an
algebraic inertia condition on $M(\eta)$. We provide the complete
proof of this correspondence, as this specific normalization is
essential for our subsequent analysis.
\end{Remark}
\begin{Corollary}\label{C:rank}
One has $\rank M(\eta)\le n+2$ for every $\ell$; for $\ell>n+2$ this
says that every $(n+3)$-minor of $M(\eta)$ vanishes. These are
universal determinantal relations among the entries of the mass
matrix.
\end{Corollary}
\begin{proof}
By Lemma~\ref{L:gram}, $-\frac12M(\eta)$ is the Gram matrix of $\ell$ vectors in the
$(n+2)$-dimensional space $\R^{1,n+1}$.
\end{proof}
For $\ell\ge n+2$ the dimensions match, at least heuristically. The symmetric
matrices of rank at most $n+2$ have codimension $\binom{\ell-n-1}{2}$ in
$\mathrm{Sym}_\ell$, and
$\binom{\ell-n-1}{2}=\dim\Hh-\big(\ell n-\tfrac12(n+1)(n+2)\big)$ is the codimension
of $\Mo_\ell(S^n)$ in $\Hh$. This is only a dimension count: $M(\eta)$ has zero
diagonal, and we do not claim that $M(\Hh)$ is transverse to the rank strata.

\begin{Corollary}\label{C:simplex}
The trivial matrix $\eta=0$ is in the image of $\Xi$ iff $\ell\le
n+2$, and then the unique such configuration is the regular simplex
with $\ell$ vertices inscribed in a round $S^{\ell-2}\subseteq S^n$.
Hence the bound $\hm_i\ge\ell-1$ is attained, simultaneously at all
ends, exactly for $\ell\le n+2$; for $\ell>n+2$ the
function $\max_i\hm_i$ attains a minimum $>\ell-1$.
\end{Corollary}

\begin{proof}
The condition $\eta=0$ gives $M=J-I$, with eigenvalues $\ell-1$
(simple) and $-1$ with multiplicity $\ell-1$; by
Theorem~\ref{T:image}(2) this is realizable iff $\ell-1\le n+1$, and
the configuration is then unique by injectivity. A regular simplex
with $\ell$ vertices has all $d_{ij}$ equal, so $\eta=0$, and its vertices
span $S^{\ell-2}$, which embeds in $S^n$ exactly when $\ell\le
n+2$. For the last statement, \eqref{eq:tropical} gives
$\max_i\hm_i\le R\Rightarrow\|\eta_{i\cdot}\|_1\le2(\ell-1)\ln R$ for
every $i$, so the sublevel sets of $\max_i\hm_i$ are preimages under
$\Xi$ of compact subsets of $\Hh$, hence compact by
Proposition~\ref{P:proper}; the minimum is therefore attained, and it
exceeds $\ell-1$ because $\eta=0$ is not in the image.
\end{proof}

\begin{Remark}\label{R:lorentzMDS}
Theorem~\ref{T:image} is the analogue, for configurations on the sphere modulo the
conformal group, of the classical embeddability criterion for distance matrices. Let
$J=I-\frac1\ell\mathbf 1\mathbf 1^{\top}$. A symmetric matrix $D$ with zero diagonal
is the squared-distance matrix of $\ell$ points in some Euclidean space if and only if
$B:=-\frac12JDJ$ is positive semidefinite \cite{Sch35,YH}, equivalently if and only if
$D$ is conditionally negative definite; such a $D\ne0$ has exactly one positive
eigenvalue, and the points are recovered from $B$, their centred Gram matrix, by its
spectral factorization (classical multidimensional scaling \cite{To}). The passage
$L\mapsto\eta$ is the analogous operation on log-distances: both $S\mapsto JSJ$ and
$\Pi$ kill exactly the additive matrices $(a_i+a_j)$ (Lemma~\ref{L:proj}). Here the role
of the squared-distance matrix is played by $M(\eta)=(e^{-2\eta_{ij}/\kp})$, and the
criterion ``exactly one positive and at most $n+1$ negative eigenvalues'' is the
statement that $M(\eta)$ is a Gram matrix of future null vectors in $\R^{1,n+1}$, the
sphere $S^n$ being the space of null rays. In this sense the mass map is a conformal,
or Lorentzian, version of multidimensional scaling, with the points recovered from
$\eta$ by Lemma~\ref{L:gram} as they are recovered from $B$ by its spectral
factorization. The relation between the two is direct: $M(\eta)$ is the Euclidean
squared-distance matrix $D=(d_{ij}^2)$ up to diagonal congruence
(Lemma~\ref{L:gram}), which is what the proof of Proposition~\ref{P:contract} exploits.
\end{Remark}

\section{The polygon model}\label{s:polygons}
This section proves Theorem~\ref{thm:mainB}. The existence of a closed
polygon with prescribed Gauss map is due to Kapovich--Leeb--Millson
\cite[Section~5.3]{KLM} and \cite{KLM2} for all symmetric spaces of noncompact
type. Balser~\cite{Ba} gives a self-contained proof for Hadamard
spaces via the fixed-point method we follow, and the map $\mathcal T$ from
Section~\ref{ss:bij} is his map $\Phi_c$ with all weights equal to
$s$. That the Gauss image of a closed polygon is nice semistable is
proved in~\cite[Section~5]{KLM}; see also~\cite{KLM2}. For
$n=2$, the entire statement, together with the symplectic structure on
the space of polygons, was proved by
Kapovich--Millson--Treloar~\cite{KMT}; the Euclidean analogue is given
in~\cite{KM96,Kl}. Our contribution is the uniqueness of the polygon in
the stable case: by the explicit step formula of Lemma~\ref{L:step}, the
map $\mathcal T$ strictly decreases distances (Lemma~\ref{L:contract}). We also prove a formula for the masses in
terms of the diagonals of a polygon (Proposition~\ref{P:massespoly}).

\subsection{The Gauss map}
Realize $\HH^{n+1}=\{x\in\R^{1,n+1}:\langle x,x\rangle=-1,\ x_0>0\}$, with
$\cosh d(x,y)=-\langle x,y\rangle$ and $\partial_\infty\HH^{n+1}=S^n$ as in
Section~\ref{s:lightcone}. Fix $s>0$ and $\lambda:=e^{-s}$.

\begin{Definition}
An \emph{$s$-equilateral $\ell$-gon} is
$x=(x_1,\dots,x_\ell)\in(\HH^{n+1})^\ell$ with $d(x_i,x_{i+1})=s$
cyclically. We write $\Polys\subset(\HH^{n+1})^\ell$ for the space of
these, and $\Po_\ell(s):=\Polys/\Isom(\HH^{n+1})$.
\end{Definition}
The moduli space $\Po_\ell(s)$ is compact: indeed, normalizing $x_1=o$
forces all vertices to lie in the ball $\overline B(o,\ell s)$, and
the residual group $\mathrm{O}(n+1)$ is compact.

\begin{Lemma}\label{L:w}
Put $w_i:=x_{i+1}-\lambda x_i$. Then $w_i$ is future null, $-\langle
w_i,x_i\rangle=\sinh s$, and $\R_{>0}w_i$ is the ideal point $p_i$
obtained by extending the side $\overrightarrow{x_ix_{i+1}}$ beyond
$x_{i+1}$, this is the Gauss image of the $i$-th side.
\end{Lemma}

\begin{proof}
We have $\langle w_i,w_i\rangle=-1+2\lambda\cosh s-\lambda^2=0$ since
  $2\lambda\cosh s=1+\lambda^2$, and $-\langle w_i,x_i\rangle=\cosh
  s-\lambda=\sinh s>0$.  If $u_i$ is the unit tangent at $x_i$ towards
  $x_{i+1}$ then $x_{i+1}=\cosh s\,x_i+\sinh s\,u_i$ and the ideal
  endpoint of $u\mapsto\cosh u\,x_i+\sinh u\,u_i$ is the ray of
  $x_i+u_i=w_i/\sinh s$.
\end{proof}

\begin{Lemma}\label{L:recon}
For any vectors $w_i$ the system $x_{i+1}=\lambda x_i+w_i$ has the
unique solution
$x_i=(1-\lambda^\ell)^{-1}\sum_{k=1}^{\ell}\lambda^{k-1}w_{i-k}$. If
the $w_i$ are null and $\langle x_i,x_i\rangle=-1$ for all $i$, then,
after replacing $(x_i,w_i)$ by $(-x_i,-w_i)$ if necessary, $x$ is an
$s$-equilateral $\ell$-gon in $\HH^{n+1}$ with sides $w_i$, and the
$w_i$ are future null.
\end{Lemma}

\begin{proof}
The cyclic system has matrix $I-\lambda C$, invertible since
$\lambda<1$, with the stated inverse. If $\langle w_i,w_i\rangle=0$
and $\langle x_i,x_i\rangle=-1$ then $0=-1-\lambda^2-2\lambda\langle
x_i,x_{i+1}\rangle$, so $-\langle x_i,x_{i+1}\rangle=\cosh s>0$; hence
all $x_i$ lie in one sheet of the hyperboloid, which the sign change
puts in $\HH^{n+1}$, and then $-\langle w_i,x_i\rangle=\sinh s>0$
shows that $w_i$ is future.
\end{proof}

Let $G_{ij}=\langle w_i,w_j\rangle$; the closing condition on our
polygon translates to an \emph{affine} one:
\begin{equation}\label{eq:closing}
\langle x_i,x_i\rangle=\frac{1}{(1-\lambda^\ell)^2}
\sum_{k,m=1}^{\ell}\lambda^{k+m-2}G_{i-k,\,i-m}=-1 .
\end{equation}

\begin{Theorem}\label{T:model}
The correspondence $x\mapsto G$ yields a homeomorphism of
$\Po_\ell(s)$ onto
\[
\mathcal G_{n,\ell}(s)=\left\{G\in\mathrm{Sym}_\ell:\
\begin{array}{l}
G_{ii}=0,\ G_{ij}\le0,\\
G\ \text{has}\le n+1\ \text{positive},\ \le1\ \text{negative eigenvalues},\\
\text{the }\ell\ \text{affine equations }\eqref{eq:closing}
\end{array}\right\}.
\]
Moreover $\rank G= 1+ \dim\{\mbox{the smallest totally geodesic subspace containing the
polygon $x$}\}$, in particular $\rank G\le n+2$, and $\rank G=2$ iff
the polygon lies on a geodesic, (which for equal side lengths requires
$\ell$ to be even).
\end{Theorem}

\begin{proof}
The map is well defined and continuous, and by Witt's theorem a
family of vectors in the nondegenerate $\R^{1,n+1}$ is determined by
its Gram matrix up to $\mathrm{O}(1,n+1)$; the future-cone condition,
guaranteed by $G_{ij}\le0$, reduces this to $\Isom(\HH^{n+1})$, so
the map is injective. Conversely the inertia condition realizes $G$ as
a Gram matrix of null vectors, and $G_{ij}\le0$ places them in one
closed cone: two null vectors with $\langle w_i,w_j\rangle=0$ are
proportional, a negatively proportional pair $w_j=-cw_i$ would force
$G_{ik}=-c^{-1}G_{jk}$ to vanish for all $k$, i.e.\ a zero row, and a
zero row is excluded by \eqref{eq:closing}: if the $i$-th row of $G$
vanishes, then $\langle
x_i,w_i\rangle=(1-\lambda^\ell)^{-1}\sum_k\lambda^{k-1}G_{i-k,i}=0$
and $\langle w_i,w_i\rangle=0$, whence $\langle
x_{i+1},x_{i+1}\rangle=\lambda^2\langle x_i,x_i\rangle=-\lambda^2\ne-1$.
Then Lemma~\ref{L:recon} produces the polygon from
\eqref{eq:closing}, so the map is a bijection; a continuous bijection
from the compact space $\Po_\ell(s)$ to the Hausdorff space
$\mathrm{Sym}_\ell$ is a homeomorphism onto its image. The span of the
$w_i$
equals the span of the $x_i$; if it is a timelike $2$-plane the
polygon is on a geodesic and its ideal points lie in a two-point set
$\{q_\pm\}$, so the signed displacement vanishes only if half the
sides point each way.
\end{proof}

\subsection{Bijectivity of the Gauss map}\label{ss:bij}
By Theorem~\ref{T:model}, the fibre of the Gauss map
$\Gamma(x)=([w_1],\dots,[w_\ell])$ over a configuration $(p_i)$ is the
set of scalings $t\in\R^\ell_{>0}$ solving the $\ell$ equations
\eqref{eq:closing} for $G_{ij}=t_it_j\langle\hat p_i,\hat p_j\rangle$,
a diagonal scaling problem of Sinkhorn type. We shall prove
bijectivity geometrically, but first we fix the target of the Gauss
map.

\subsubsection*{Stability and the polystable quotient}
Let $X:=(S^n)^\ell$. For $\pl\in X$ and $\zeta\in S^n$ write
$k_\zeta:=\#\{i:p_i=\zeta\}$. We call the equally weighted
configuration $\pl$ \emph{stable} if $k_\zeta<\ell/2$ for every
$\zeta$, \emph{semistable} if $k_\zeta\le\ell/2$ for every $\zeta$,
\emph{strictly semistable} if it is semistable but not stable,
i.e.\ $k_\zeta=\ell/2$ for some $\zeta$, and \emph{polystable} (or
nice semistable) if it is stable or, when $\ell$ is even, it is
supported at exactly two points with multiplicity $\ell/2$ at each.
Let $\Xst\subset\Xss\subset X$ be the stable and semistable loci;
they are $\Mob(S^n)$-invariant, and $\Conf_\ell(S^n)\subset\Xst$
since $k_\zeta\le1<\ell/2$. For $\ell$ even and a splitting
$[\ell]:=\{1,\dots,\ell\}=A\sqcup B$ with $|A|=|B|=\ell/2$, let $\pAB$
denote any configuration with the points of $A$ at one point of $S^n$
and those of $B$ at another. Since $\Mob(S^n)$ is transitive on pairs
of distinct points, these configurations form a single orbit
$X_{\{A,B\}}$; they are the non-stable polystable configurations. We
need one more moduli space:
\[
\overline{\Mo}^{\,\mathrm{ps}}_\ell(S^n)
:=\{\text{polystable equally weighted configurations}\}/\Mob(S^n),
\]
with the topology of the semistable quotient of \cite{KLM}: the
quotient $q\colon\Xss\to\MopsL$ of the semistable locus in which each
orbit is identified with the unique closed orbit in its closure, which
is Hausdorff by \cite[Section~4]{KLM}. Thus $\MopsL$ is the closed-orbit
model of the semistable quotient, not the raw orbit space of all
semistable configurations. That the closed orbits in $\Xss$ are
exactly the polystable ones is Lemma~\ref{L:orbits} below; it relies
on the convergence property of M\"obius groups, which we recall.

\begin{Lemma}[divergent sequences in $\Mob(S^n)$]\label{L:northsouth}
Let $(g_m)$ be a sequence in $\Mob(S^n)$ with no convergent
subsequence. Then, after passing to a subsequence, there are points
$\zeta_+,\zeta_-\in S^n$, not necessarily distinct, such that
$g_m\to\zeta_+$ uniformly on compact subsets of
$S^n\setminus\{\zeta_-\}$.
\end{Lemma}
This is the convergence property of M\"obius groups, see
\cite{GM}; it also follows directly from the Cartan decomposition
$g=ka_tk'$ of $\Isom(\HH^{n+1})$, $k,k'\in\mathrm{O}(n+1)$, $a_t$ a
translation of length $t=d(o,go)$ along a fixed geodesic.

\begin{Lemma}[closed orbits in the semistable locus]\label{L:orbits}
Let $\ell\ge3$.
\begin{enumerate}[label=\textup{(\arabic*)},leftmargin=*]
\item If $\pl$ is stable, then the orbit $\Mob(S^n)\,\pl$ is closed in $\Xss$.
\item If $\ell$ is even, each orbit $X_{\{A,B\}}$ is closed in $\Xss$.
\item Let $\pl$ be strictly semistable, let $A:=\{i:p_i=\zeta\}$ be
  its cluster of $\ell/2$ points, and $B:=[\ell]\setminus A$. Then
  $X_{\{A,B\}}$ is the unique closed orbit contained in
  \[
  \overline{\Mob(S^n)\,\pl}\cap\Xss.
  \]
  In particular $\Mob(S^n)\,\pl$
  is closed in $\Xss$ only if $\pl\in X_{\{A,B\}}$.
\end{enumerate}
Thus the closed $\Mob(S^n)$-orbits in $\Xss$ are exactly the orbits of
polystable configurations, and the quotient map satisfies $q(\pl)=\Mob(S^n)\,\pl$ for
$\pl\in\Xst$ and $q(\pl)=q(\pAB)$ for $\pl$ as in \textup{(3)}.
\end{Lemma}

\begin{proof}
Since $X$ is metrizable, closedness can be tested on sequences.

(1) Let $\pl$ be stable and $g_m\pl\to\ql\in\Xss$. If some subsequence
$g_{m_j}$ converges to $g\in\Mob(S^n)$, then $\ql=g\pl$. Otherwise
Lemma~\ref{L:northsouth} applies: along a subsequence,
$g_mp_i\to\zeta_+$ for every $i$ with $p_i\ne\zeta_-$, so
$q_i=\zeta_+$ for all such $i$. There are
$\ell-k_{\zeta_-}(\pl)>\ell/2$ of them by stability, so
$k_{\zeta_+}(\ql)>\ell/2$ and $\ql$ is unstable, a contradiction.

(2) Let $\pl=\pAB$ with $A$ at $\zeta$ and $B$ at $\zeta'\ne\zeta$,
and let $g_m\pl\to\ql\in\Xss$. If $(g_m)$ has a convergent
subsequence, $\ql$ lies in the orbit as in (1). Otherwise pass to a
subsequence as in Lemma~\ref{L:northsouth} and such that $g_m\zeta$
and $g_m\zeta'$ converge. If $\zeta_-\notin\{\zeta,\zeta'\}$, all
$\ell$ points of $\ql$ equal $\zeta_+$, which is unstable. If
$\zeta_-=\zeta$, say, then in $\ql$ the points of $B$ sit at $\zeta_+$
and those of $A$ at $\zeta'':=\lim g_m\zeta$; semistability forces
$\zeta''\ne\zeta_+$, so $\ql\in X_{\{A,B\}}$.

(3) Normalize $\zeta=0$ in the upper half-space model and let
$h_t\in\Mob(S^n)$ act on $\R^n\cup\{\infty\}$ by $x\mapsto
e^{t}x$. Then $h_tp_i=0$ for $i\in A$, while $h_tp_j\to\infty$ for
$j\in B$ because $p_j\ne0$. Hence $h_t\pl\to\pAB$ with $A$ at $0$ and
$B$ at $\infty$, a configuration in $\Xss$; so $X_{\{A,B\}}$ lies in
$\overline{\Mob(S^n)\,\pl}\cap\Xss$, and it is closed by (2). For
uniqueness let $\Mob(S^n)\,\ql$ be a closed orbit in
\[
\overline{\Mob(S^n)\,\pl}\cap\Xss.
\]
The condition that the points of $A$ coincide is closed and
$\Mob(S^n)$-invariant, so it holds for $\ql$, and $\ql$ is not
stable. If $\ql\notin X_{\{A,B\}}$, then $\ql$ is strictly semistable
with $A$ as its cluster of $\ell/2$ points, and what we just proved,
applied to $\ql$, shows that $\overline{\Mob(S^n)\,\ql}\cap\Xss$
contains the orbit $X_{\{A,B\}}$, which is different from
$\Mob(S^n)\,\ql$ because $X_{\{A,B\}}$ is a single orbit and
$\ql\notin X_{\{A,B\}}$; this contradicts the closedness of
$\Mob(S^n)\,\ql$. Hence $\ql\in X_{\{A,B\}}$.

For the consequence, a closed orbit is stable or strictly semistable, and in the
latter case (3) shows that it is some $X_{\{A,B\}}$; conversely these orbits are
closed by (1) and (2). The description of $q$ is then the definition of $\MopsL$ as
the closed-orbit quotient.
\end{proof}

In the language of \cite{KLM}, Lemma~\ref{L:orbits} verifies directly, for equal
weights on $S^n$, that the nice semistable configurations are those with closed
orbit.

\subsubsection*{The fixed-point method}
For $p\in S^n$ let $\mathcal T_p(x)$ be the point at distance $s$ from $x$
on the ray $[x,p)$, and $\mathcal T:=\mathcal T_{p_\ell}\circ\cdots\circ\mathcal T_{p_1}$.

\begin{Lemma}\label{L:fixed}
The correspondence $x_1\mapsto(x_1,\mathcal T_{p_1}(x_1),\dots)$ is a
bijection between fixed points of $\mathcal T$ and $s$-equilateral
$\ell$-gons with Gauss image $(p_1,\dots,p_\ell)$.
\end{Lemma}
We need just one computation. Normalize $\HH^{n+1}$ as the upper
half-space $\{(z,t):z\in\R^n,\ t>0\}$ with $ds^2=(|dz|^2+dt^2)/t^2$
and ideal boundary $\R^n\cup\{\infty\}$.

\begin{Lemma}
  \label{L:step}
Let $p\in\R^n$, $x=(z,t)$, $a=|z-p|$, $\sigma=a/t$, and write
$\mathcal T_p(x)=(z',t')$, $a'=|z'-p|$. Then $z'$ lies on the segment
$[z,p]$, and
\[
t'=\mu(\sigma)\,t,\qquad a'=\nu(\sigma)\,a,\qquad
\mu(\sigma)=\frac{(1+\sigma^2)e^{s}}{e^{2s}+\sigma^2},\quad
\nu(\sigma)=\frac{1+\sigma^2}{e^{2s}+\sigma^2},
\]
so that $a'/t'=e^{-s}\sigma$. Here $\mu$ increases from
$\mu(0)=e^{-s}$ to $\mu(\infty)=e^{s}$, with $\mu(\sigma)=1$ iff
$\sigma=e^{s/2}$, and $\nu$ increases from $e^{-2s}$ to $1$. For
$p=\infty$ one has $\mathcal T_\infty(z,t)=(z,e^{s}t)$.
\end{Lemma}

\begin{proof}
For $a=0$, and for $p=\infty$, the geodesic is vertical and the length element is
$dt/t$. Let $a>0$. The geodesic through $x$ with ideal endpoint $p$ lies in the
vertical $2$-plane containing $z$ and $p$; in coordinates $(u,t)$ of that plane, with
$p$ at $u=0$ and $x$ at $u=a$, it is the half-circle through the origin centred at
$(c,0)$ with $c=r=(a^2+t^2)/2a$. Parametrizing by
$(u,t)=(c+r\cos\theta,r\sin\theta)$ the length element is $d\theta/\sin\theta$, so
$\ln\tan(\theta/2)$ is an arclength parameter increasing towards $\theta=\pi$, which
is the endpoint $p$; moving a distance $s$ towards $p$ therefore multiplies
$\tan(\theta/2)$ by $e^{s}$. At $x$ one computes $\tan(\theta_0/2)=t/a$, so
$T:=\tan(\theta_1/2)=e^{s}t/a$, and with $\sin\theta=2T/(1+T^2)$,
$1+\cos\theta=2/(1+T^2)$,
\[
t'=r\sin\theta_1=\frac{(a^2+t^2)te^{s}}{a^2+t^2e^{2s}},\qquad
a'=r(1+\cos\theta_1)=\frac{a(a^2+t^2)}{a^2+t^2e^{2s}} ,
\]
which are the stated formulas after dividing by $t^2$. Since
$\theta\mapsto c+r\cos\theta$ decreases monotonically from $a$ to $0$ on
$[\theta_0,\pi]$, the horizontal coordinate moves from $z$ to $p$ without
overshooting. The monotonicity of $\mu$ and $\nu$ is read off
$\mu(\sigma)=e^{s}\big(1-\frac{e^{2s}-1}{e^{2s}+\sigma^2}\big)$ and $\nu=e^{-s}\mu$.
\end{proof}

The exact law $a'/t'=e^{-s}\sigma$ has a coordinate-free meaning: since
$a/t=\sinh d(x,\gamma)$ for $\gamma$ the geodesic from $p$ to $\infty$, moving a
distance $s$ towards $p$ contracts by exactly $e^{-s}$ the hyperbolic sine of the
distance to any geodesic ending at $p$. It is the hyperbolic sine of the
distance, not the distance itself, that is contracted by $e^{-s}$.

\begin{Lemma}\label{L:contract}
The function $\mathcal T_p$ is nonexpansive, and $d(\mathcal T_px,\mathcal T_py)=d(x,y)$
if and only if $x$ and $y$ lie on one geodesic with ideal endpoint
$p$. Hence $\mathcal T$ is nonexpansive; and if the $p_i$ take at least
three distinct values then $d(\mathcal T x,\mathcal T y)<d(x,y)$ for $x\ne y$, so
$\mathcal T$ has at most one fixed point.
\end{Lemma}

\begin{proof}
Normalize $p=\infty$, so $\mathcal T_p(z,t)=(z,e^{s}t)$. From
$\cosh d((z,t),(w,\tau))=1+\frac{|z-w|^2+(t-\tau)^2}{2t\tau}$,
\[
\cosh d\big(\mathcal T_p(z,t),\mathcal T_p(w,\tau)\big)
=1+\frac{e^{-2s}|z-w|^2+(t-\tau)^2}{2t\tau},
\]
which is $\le\cosh d((z,t),(w,\tau))$, with equality iff $z=w$, i.e.\ iff the points
lie on a common geodesic ending at $p$ --- along which $\mathcal T_p$ is the restriction of
a hyperbolic translation, hence an isometry. If $x\ne y$ and every factor of $\mathcal T$
preserved the distance of the running pair, the geodesic $g$ through $x$ and $y$ would
be preserved by each $\mathcal T_{p_i}$ and would have every $p_i$ among its two ideal
endpoints, which is impossible with three distinct values.
\end{proof}

\begin{Theorem}[Theorem~\ref{thm:mainB}]\label{T:gauss}
Let $p_1,\dots,p_\ell\in S^n$, $\ell\ge3$.
\begin{enumerate}[label=\textup{(\arabic*)},leftmargin=*]
\item If the $p_i$ take at least three distinct values, there is at most one
$s$-equilateral $\ell$-gon with Gauss image $(p_i)$, up to isometry.
\item If $k_\zeta>\ell/2$ for some $\zeta$, there is none.
\item If $k_\zeta<\ell/2$ for every $\zeta$, there is exactly one, up to isometry.
\item If the configuration is supported at exactly two points, then a polygon
exists exactly when the two multiplicities are equal; in that case all fixed points
lie on the joining geodesic and determine a single point of $\Po_\ell(s)$.
\item If $k_\zeta=\ell/2$ for some $\zeta$ and the $p_i$ take at least three
distinct values, there is none.
\end{enumerate}
Consequently the Gauss map induces a homeomorphism
\[
\Gamma:\Po_\ell(s)\xrightarrow{\ \cong\ }
\overline{\Mo}^{\,\mathrm{ps}}_\ell(S^n).
\]
On the open locus of polygons with pairwise distinct ideal points it
restricts to a homeomorphism onto $\Mo_\ell(S^n)$
(Corollary~\ref{C:compact}).
\end{Theorem}

Note that the hypothesis of (1) is implied by that of (3): two distinct
values would give $k_{\zeta_1}+k_{\zeta_2}=\ell$ with both terms $<\ell/2$. It cannot
be dropped, by Remark~\ref{R:degenerate}.

\begin{proof}
(1) is Lemmas~\ref{L:fixed} and \ref{L:contract}.

(2) Normalize $\zeta=\infty$ and put $k=k_\infty>\ell/2$. Each of the $k$ factors with
$p_i=\infty$ multiplies $t$ by exactly $e^{s}$, and each of the other $\ell-k$ factors
multiplies it by $\mu(\sigma)\ge\mu(0)=e^{-s}$; hence
$t(\mathcal T x)\ge e^{(2k-\ell)s}t(x)>t(x)$ for every $x$, and a fixed point would have
$t(\mathcal T x)=t(x)$.

(5) Normalize $\zeta=\infty$. Now the $\ell/2$ factors with
$p_i=\infty$ multiply $t$ by $e^{s}$ and the others by $\mu(\sigma)\ge
e^{-s}$, with equality iff $\sigma=0$, i.e.\ iff the running point has
horizontal coordinate $p_i$. Hence $t(\mathcal T x)\ge t(x)$, with
equality iff at every step towards a finite $p_i$ the running
horizontal coordinate equals $p_i$. But neither $\mathcal T_\infty$ nor
$\mathcal T_p$ at $\sigma=0$ moves the horizontal coordinate, so equality
forces all finite $p_i$ to coincide: only two distinct values.

(3) The $p_i$ take $d\ge3$ distinct values $q_1,\dots,q_d$ with
multiplicities $k_\alpha<\ell/2$. Normalize so that
$\infty\notin\{p_i\}$, and put $P=\max_i|p_i|$ and
$\delta_0=\min_{q_\alpha\ne q_\beta}|q_\alpha-q_\beta|>0$. Let
$w_0=(z_0,t_0)$ be arbitrary, $w_{r+1}=\mathcal T_{p_{\pi(r)}}(w_r)$ with
$\pi(r)\equiv r+1\pmod\ell$, so that $w_{m\ell}=\mathcal T^mw_0$; any $\ell$
consecutive steps use each $p_i$ once. We show the orbit is bounded,
in three steps.

\emph{The horizontal coordinate never escapes.} By Lemma~\ref{L:step}, each
$z_{r+1}$ lies on the segment $[z_r,p_{\pi(r)}]$. Hence
$|z_{r+1}|\le\max(|z_r|,P)$ and $|z_r|\le Z:=\max(|z_0|,P)$ for all $r$.
Consequently, with $C:=Z+P$,
\begin{equation}\label{eq:sigmaC}
\begin{array}{c}
  \sigma_r:=\frac{|z_r-p_{\pi(r)}|}{t_r}\ \le\ \frac{C}{t_r} .
\end{array}
\end{equation}
\emph{The parameter $t$ is bounded above.} If $t_r>Ce^{-s/2}$ then
$\sigma_r<e^{s/2}$ by \eqref{eq:sigmaC}, so $\mu(\sigma_r)<1$ and
$t_{r+1}<t_r$; and if $t_r\le Ce^{-s/2}$ then $t_{r+1}\le e^{s}t_r\le
Ce^{s/2}$. Hence $t_r\le\max(t_0,Ce^{s/2})$ for all $r$.

\emph{The parameter $t$ is bounded below; this is where stability enters.} Fix
\[
\epsilon:=e^{-\ell
  s}\,\delta_0\,\min\Big\{\tfrac18,\ \sqrt{\tfrac{1-e^{-s}}{128\,\ell\,(e^{2s}-1)}}\Big\},
\]
and suppose $t_r\le\epsilon$. Every step multiplies $t$ by at most $e^{s}$, so
$t_{r+j}\le\tau:=e^{\ell s}\epsilon\le\delta_0/8$ for $0\le j\le\ell$, and
we obtain
\[
\begin{array}{c}
\tau^2\ \le\ \frac{\delta_0^2(1-e^{-s})}{128\,\ell\,(e^{2s}-1)}\ \le\ \frac{\delta_0^2}{64\,\ell\,(e^{2s}-1)} .
\end{array}
\]
Let $q_*$ be the value nearest to $z_r$ and $k$ its multiplicity if
$|z_r-q_*|<\delta_0/4$, and $k=0$ otherwise. By Lemma~\ref{L:step} a
step towards a target at horizontal distance $a$ moves the horizontal
coordinate along the segment to the target by
$(1-\nu(\sigma))a\le(e^{2s}-1)t^2/a$. We claim that throughout the
pass every target other than $q_*$ (every target, if $k=0$) stays at
distance $\ge\delta_0/8$. Indeed, as long as this holds the drift of
each step towards such a target is at most
$8(e^{2s}-1)\tau^2/\delta_0$, so the total drift over the pass is at
most $8\ell(e^{2s}-1)\tau^2/\delta_0\le\delta_0/8$; steps towards
$q_*$ do not increase $|z-q_*|$, so $|z-q_*|\le3\delta_0/8$ throughout
and the other targets stay at distance $\ge5\delta_0/8$ if $k>0$,
while if $k=0$ all targets, initially at distance $\ge\delta_0/4$,
stay at distance $\ge\delta_0/8$. The claim follows by induction over
the pass. Now the $k$ steps towards $q_*$ contribute at least
$\mu(0)^k=e^{-ks}$, while each of the other $\ell-k$ steps has
$\sigma\ge\delta_0/8\tau$ and so contributes
\[
\begin{array}{c}
\mu(\sigma)\ \ge\ e^{s}\Big(1-\frac{64(e^{2s}-1)\tau^2}{\delta_0^2}\Big)
\ \ge\ e^{s}\Big(1-\frac{1-e^{-s}}{2\ell}\Big).
\end{array}
\]
Hence, using $(1-u)^\ell\ge1-\ell u$ and $\ell-2k\ge1$,
\[
\begin{array}{c}
\frac{t_{r+\ell}}{t_r}\ \ge\ e^{(\ell-2k)s}\Big(1-\frac{1-e^{-s}}{2}\Big)
\ \ge\ \frac{1+e^{s}}{2}\ >\ 1 .
\end{array}
\]
A first-descent argument now gives $t_r\ge e^{-\ell s}\min(t_0,\epsilon)$ for all $r$:
if $r$ is minimal with $t_r$ below that bound then $r\ge\ell$, and $t_{r-\ell}$ either
is $\le\epsilon$, so that the pass estimate applies, or exceeds $\epsilon$, in which
case $t_r\ge e^{-\ell s}t_{r-\ell}$; either way a contradiction.

The orbit therefore lies in a bounded subset of $\HH^{n+1}$. Let
$x_m:=\mathcal T^m(w_0)$ and
\[
r(y):=\limsup_{m\to\infty}d(y,x_m).
\]
A bounded sequence in the complete $\mathrm{CAT}(0)$ space $\HH^{n+1}$ has a unique
asymptotic center, i.e.\ a unique minimizer $c$ of $r$. Since $\mathcal T$ is nonexpansive,
\[
r(\mathcal T c)=\limsup_m d(\mathcal T c,x_m)
\le\limsup_m d(c,x_{m-1})=r(c).
\]
Thus $\mathcal T c$ is also an asymptotic center, and uniqueness gives $\mathcal T c=c$.
This is the fixed-point argument of \cite[Section~5]{KLM} and \cite{Ba}.
Uniqueness of the resulting polygon in the stable case is (1).

For (4), let the two ideal points be $\zeta_\pm$ with multiplicities
$k_\pm$; relabelling, we may assume $k_+\ge k_-$. Normalize
$\zeta_+=\infty$ and $\zeta_-=0$, so that the joining geodesic
$\gamma$ is the vertical line $z=0$. By Lemma~\ref{L:step},
$t(\mathcal T x)\ge e^{(k_+-k_-)s}\,t(x)$ for every $x$, with equality if
and only if the running horizontal coordinate is $0$ at each of the
$k_-$ steps towards $0$; since neither $\mathcal T_\infty$ nor $\mathcal T_0$ at
$\sigma=0$ moves the horizontal coordinate, equality holds if and only
if $x\in\gamma$. Hence for $k_+>k_-$ there is no fixed point, while
for $k_+=k_-=\ell/2$ every fixed point lies on $\gamma$, on which
$\mathcal T$ acts as the translation of signed length $(k_+-k_-)s=0$: every
point of $\gamma$ is fixed. Translations along $\gamma$ act
transitively on these fixed points and stabilize the two ideal points,
so they determine one point of $\Po_\ell(s)$.

Parts (2), (4) and (5) say that the Gauss image of an $s$-equilateral polygon is
polystable; (3) and (4) that every polystable configuration is a Gauss image; (1) and
(4) that the polygon is unique up to isometry. Hence $\Gamma$ is a
bijection between $\Po_\ell(s)$ and the set of polystable orbits,
which by Lemma~\ref{L:orbits} is the set underlying $\MopsL$. The
Gauss map $\Polys\to\Xss$, $x\mapsto([w_1],\dots,[w_\ell])$, is
continuous and $\Isom(\HH^{n+1})$-equivariant, so its composite with
$q$ descends to a continuous map $\Gamma\colon\Po_\ell(s)\to\MopsL$;
since $\Po_\ell(s)$ is compact and $\MopsL$ is Hausdorff \cite{KLM},
the bijection is a homeomorphism. The last assertion of the theorem is
part of Corollary~\ref{C:compact} below.
\end{proof}

\begin{Corollary}[the polystable quotient as a compactification]
  \label{C:compact}
Let $\ell\ge3$.
\begin{enumerate}[label=\textup{(\arabic*)},leftmargin=*]
\item The moduli space $\MopsL$ is compact and metrizable.
\item The inclusion $\Conf_\ell(S^n)\subset\Xss$ induces a
  homeomorphism of $\Mo_\ell(S^n)$ onto an open dense subset of
  $\MopsL$; under $\Gamma$ this subset is the open locus of polygons
  with pairwise distinct ideal points. Thus $\MopsL$ is a
  compactification of $\Mo_\ell(S^n)$.
\item If $\ell$ is odd, then $\Xss=\Xst$, and $\MopsL=\Xst/\Mob(S^n)$
  is the orbit space of the stable locus.
\end{enumerate}
\end{Corollary}

\begin{proof}
We write $\Mob$ for $\Mob(S^n)$.

(1) By Theorems~\ref{T:gauss} and \ref{T:model}, $\MopsL$ is
homeomorphic to $\Po_\ell(s)$ and hence to $\Gmod\subset\mathrm{Sym}_\ell$;
the former is compact and the latter is metrizable.

(2) For $S\subset[\ell]$ let $\Delta_S:=\{\pl\in X:\ p_i=p_j\text{ for
all }i,j\in S\}$, a closed $\Mob$-invariant subset of $X$. Then
$X\setminus\Xss=\bigcup_{|S|>\ell/2}\Delta_S$ and
$X\setminus\Xst=\bigcup_{|S|\ge\ell/2}\Delta_S$, so $\Xst\subset\Xss$
are open $\Mob$-invariant subsets of $X$.

By Lemma~\ref{L:orbits}, for $\pl\in\Xst$ the closed orbit
representing $q(\pl)$ is $\Mob\,\pl$ itself, whereas for
$\pl\in\Xss\setminus\Xst$ it is an orbit $X_{\{A,B\}}$, which is not
stable. Hence, for $\pl,\pl'\in\Xst$, $q(\pl)=q(\pl')$ iff
$\Mob\,\pl=\Mob\,\pl'$, and $q(\Xst)\cap q(\Xss\setminus\Xst)=\emptyset$.
Thus $\Xst=q^{-1}(q(\Xst))$ is $q$-saturated and open, so $q(\Xst)$ is
open in $\MopsL$; and the restriction of the quotient map $q$ to the
open saturated subset $\Xst$ is a quotient map onto its image.
Therefore $q$ induces a homeomorphism of the orbit space $\Xst/\Mob$
onto the open subset $q(\Xst)$ of $\MopsL$.

Next, for any space $Y$ with a $\Mob$-action the orbit map $Y\to
Y/\Mob$ is open, the saturation $\bigcup_{g\in\Mob}gU$ of an open set
being open. If $Z\subset Y$ is open and $\Mob$-invariant, the quotient
topology of $Z/\Mob$ therefore coincides with the subspace topology of
its image in $Y/\Mob$, and that image is open. Applying this with
$Y=\Xst$ and $Z=\Conf_\ell(S^n)$ identifies
$\Mo_\ell(S^n)=\Conf_\ell(S^n)/\Mob$ with an open subspace of
$\Xst/\Mob$, hence of $\MopsL$. It is dense: $\Conf_\ell(S^n)$ is dense
in $X$, its complement being the finite union of the closed sets
$\{p_i=p_j\}$ with empty interior, hence dense in $\Xss$, and $q$ is a
continuous surjection. Finally, $\Gamma^{-1}(\Mo_\ell(S^n))$ consists
of the polygons whose Gauss image lies in $\Conf_\ell(S^n)$, i.e.\ has
pairwise distinct ideal points.

(3) If $\ell$ is odd, no subset of $[\ell]$ has exactly $\ell/2$
elements, so $\Xss=\Xst$ and $\MopsL=q(\Xst)\cong\Xst/\Mob$ by (2).
\end{proof}

\begin{Remark}\label{R:degenerate}
If the $p_i$ take exactly two values $\zeta_\pm$ with multiplicities
$k_\pm$ then $\mathcal T$ preserves the geodesic $\gamma$ joining them and
acts on it as the translation by $(k_+-k_-)s$. For $k_+=k_-=\ell/2$
every point of $\gamma$ is fixed, and the $s$-equilateral $\ell$-gons
with that Gauss image form a one-parameter family --- the degenerate
polygons of Theorem~\ref{T:model} lying on $\gamma$, permuted by the
translations along $\gamma$, which stabilize the configuration. These
are the non-stable polystable representatives. They should not be
confused with the whole strictly semistable locus: a general strictly
semistable orbit need not itself have two-point support, but its orbit
closure contains the corresponding polystable representative
(Lemma~\ref{L:orbits}(3)).
\end{Remark}

\begin{Remark}\label{R:Sl}
The homeomorphism $\Gamma:\Po_\ell(s)\xrightarrow{\ \cong\ }
\overline{\Mo}^{\,\mathrm{ps}}_\ell(S^n)$ uses the cyclic order of the
polygon, so Theorem~\ref{T:gauss} endows $\Po_\ell(s)$ with an
$S_\ell$-action of which only the dihedral subgroup is visible:
reordering the ideal points produces a genuinely different
polygon. Numerically, for $n=3$, $\ell=5$, $s=1$, three orderings of
the same five points give polygons with different multisets of vertex
distances (they differ already in their smallest element). In the
Euclidean model the action is obvious, the closing condition
$\sum\vec e_i=0$ being symmetric; for hyperbolic polygons it is a
nontrivial consequence of the theorem.
\end{Remark}

\subsubsection*{The boundary of the polystable quotient}
We describe the stratification of $\MopsL$ by coincidence patterns; it is used in
\cite{PartII} to compare $\MopsL$ with the bubble-tree compactifications.
For $\pl\in X$ the \emph{coincidence pattern} $\pi(\pl)$ is the partition of $[\ell]$
into the fibres of $i\mapsto p_i$, and for a partition $\pi$ we put
$X_\pi:=\{\pl\in X:\pi(\pl)=\pi\}$, a $\Mob(S^n)$-invariant subset. We call $\pi$
\emph{stable} if every block has fewer than $\ell/2$ elements; a configuration is
stable iff its pattern is, so $\Xst=\bigsqcup_{\pi\ \mathrm{stable}}X_\pi$, and a stable
partition has at least three blocks. We write $\Mob$ for $\Mob(S^n)$.

\begin{Proposition}[the boundary of the polystable quotient]\label{P:psboundary}
Let $\ell\ge3$.
\begin{enumerate}[label=\textup{(\arabic*)},leftmargin=*]
\item \textup{(The Gram model.)} Sending a polystable configuration to the Gram matrix
$G=(\langle w_i,w_j\rangle)$ of the sides of its $s$-polygon defines a homeomorphism
\[
\Theta_s\colon\MopsL\xrightarrow{\ \cong\ }\Gmod\subset\mathrm{Sym}_\ell ,
\]
onto the Gram model of Theorem~\ref{T:model},
and $G_{ij}=0$ if and only if $p_i=p_j$. Hence $\Theta_s$ identifies $\Mo_\ell(S^n)$
with $\{G\in\Gmod:\ G_{ij}<0\text{ for all }i\ne j\}$ and $\MopsL$ with the closure of
this set in $\mathrm{Sym}_\ell$; and $\Gmod$ is a compact semialgebraic subset of
$\mathrm{Sym}_\ell$.

\item \textup{(The boundary.)} As a set,
\[
\MopsL\;=\;\bigsqcup_{\pi\ \mathrm{stable}}X_\pi/\Mob(S^n)\ \ \sqcup
\bigsqcup_{\{A,B\}}\{\pt_{A|B}\},\qquad \pt_{A|B}:=q\big(X_{\{A,B\}}\big),
\]
the second union running over the $\tfrac12\binom{\ell}{\ell/2}$ unordered splittings
$[\ell]=A\sqcup B$ with $|A|=|B|=\ell/2$ (it is void for $\ell$ odd). Moreover:
\begin{enumerate}[label=\textup{(\alph*)},leftmargin=*]
\item each $X_\pi/\Mob(S^n)$ carries the quotient topology of $X_\pi$, is Hausdorff,
and is homeomorphic to $\Mo_{|\pi|}(S^n)$ via the configuration of the distinct values,
labelled by the blocks of $\pi$; the discrete partition contributes $\Mo_\ell(S^n)$,
and $\Xst/\Mob(S^n)=\bigsqcup_\pi X_\pi/\Mob(S^n)$ is open in $\MopsL$;
\item $q^{-1}(\pt_{A|B})$ is the set of semistable configurations in which the points
of $A$ coincide or the points of $B$ coincide; the strictly semistable locus is the
union of these finitely many fibres;
\item for every stable $\pi$,
\[
\overline{X_\pi/\Mob(S^n)}\;=\bigsqcup_{\substack{\pi'\succeq\pi\\ \pi'\ \mathrm{stable}}}
X_{\pi'}/\Mob(S^n)\ \ \sqcup\ \ \big\{\pt_{A|B}:\ A\text{ is a union of blocks of }\pi\big\},
\]
where $\pi'\succeq\pi$ means that $\pi'$ is coarser than $\pi$;
\item under $\Theta_s$, the stratum $X_\pi/\Mob(S^n)$ is the locus
$\{G\in\Gmod:\ G_{ij}=0\iff i\sim_\pi j\}$, and $\pt_{A|B}$ is the unique $G\in\Gmod$
with $G_{ij}=0$ exactly when $i$ and $j$ lie on the same side of the splitting; the
points $\pt_{A|B}$ are the points of rank $2$ of Theorem~\ref{T:model}.
\end{enumerate}
\end{enumerate}
\end{Proposition}

\begin{proof}
\emph{Step 1: the decomposition.} By Lemma~\ref{L:orbits}, for
$\pl\in\Xst$ the closed orbit representing $q(\pl)$ is $\Mob\,\pl$ itself, so
$q(\Xst)=\Xst/\Mob$ as a set. Further, $q(\Xss\setminus\Xst)$ consists of the points $\pt_{A|B}$: by
Lemma~\ref{L:orbits}(3) a strictly semistable configuration whose cluster of $\ell/2$
points is $A$ is sent to $q(\pAB)$, and $q(\pAB)=\pt_{A|B}$ is a single point
because $X_{\{A,B\}}$ is one orbit. The points $\pt_{A|B}$ are pairwise distinct, their
closed orbits having different coincidence patterns, and they are distinct from the
points of $\Xst/\Mob$. Together with $\Xst=\bigsqcup_{\pi\ \mathrm{stable}}X_\pi$ this
gives the decomposition in (2). For (2b): if $\pl$ is semistable and the points of $A$
coincide at $\zeta$, then $k_\zeta=\ell/2$ exactly, no point of $B$ being allowed at
$\zeta$, so $\pl$ is strictly semistable with cluster $A$ and $q(\pl)=\pt_{A|B}$;
conversely, if $q(\pl)=\pt_{A|B}$ then $\pl\notin\Xst$, so $\pl$ is strictly
semistable with some cluster $A'$ of $\ell/2$ points, and $q(\pl)=\pt_{A'|[\ell]\setminus A'}$
forces $\{A',[\ell]\setminus A'\}=\{A,B\}$. Since a strictly semistable configuration
has a unique cluster of $\ell/2$ points unless it lies in some $X_{\{A,B\}}$, the
strictly semistable locus is the union of the fibres $q^{-1}(\pt_{A|B})$.

\emph{Step 2: the strata.} As in the proof of Corollary~\ref{C:compact}(2), the
orbit map $Y\to Y/\Mob$ of a $\Mob$-space is open; for a $\Mob$-invariant (not
necessarily open) $Z\subset Y$ the image of $V\cap Z$, $V$ open in $Y$, is the
intersection of the open image of $V$ with the image of $Z$, so the quotient topology
of $Z/\Mob$ is the subspace topology of its image. With $Y=\Xst$ and $Z=X_\pi$,
$\pi$ stable, $X_\pi/\Mob$ is a subspace of $\Xst/\Mob\cong q(\Xst)\subset\MopsL$,
hence Hausdorff. Ordering the blocks of $\pi$ by their least elements, the map
$X_\pi\to\Conf_{|\pi|}(S^n)$ recording the common value of each block is a
$\Mob$-equivariant homeomorphism, whence
$X_\pi/\Mob\cong\Conf_{|\pi|}(S^n)/\Mob=\Mo_{|\pi|}(S^n)$. For the discrete partition,
$X_\pi=\Conf_\ell(S^n)$ is open in $\Xst$, so $\Mo_\ell(S^n)$ is open in $\Xst/\Mob$ and
hence in $\MopsL$. This proves (2a).

\emph{Step 3: the Gram model.} $\Theta_s$ is the composite of $\Gamma^{-1}$, the inverse of
the Gauss map of Theorem~\ref{T:gauss}, with the homeomorphism of
Theorem~\ref{T:model}, hence a homeomorphism. If $x$ is an $s$-polygon
with sides $w_i=t_i\hat p_i$, $t_i>0$, then by \eqref{eq:chordal}
\[
G_{ij}=\langle w_i,w_j\rangle=-\tfrac12\,t_it_j\,|p_i-p_j|^2 ,
\]
which vanishes iff $p_i=p_j$. So the zero pattern of $\Theta_s(q(\pl))$ is the
coincidence pattern of the Gauss image: $\Mo_\ell(S^n)$ corresponds to the matrices
with all $G_{ij}<0$, the stratum $X_\pi/\Mob$ to the zero pattern $\pi$, and
$\pt_{A|B}$ to the zero pattern $\{A,B\}$, where Theorem~\ref{T:gauss}(4) says that
exactly one class of polygons has a Gauss image of this pattern. For such a polygon
the sides span the timelike $2$-plane of two null directions, so $\rank G=2$, and
conversely the rank-$2$ points of $\Gmod$ are the polygons on a geodesic
(Theorem~\ref{T:model}), whose Gauss image has two-point support. This proves (2d) and
the first two assertions of (1). Since $\Gmod$ is compact and contains
$\Theta_s(\Mo_\ell(S^n))$ as a dense subset, it is the closure of that subset in
$\mathrm{Sym}_\ell$. Finally, $\Gmod$ is semialgebraic: it is cut out by the polynomial
conditions $G_{ii}=0$, $G_{ij}\le0$, the affine closing equations of \eqref{eq:closing}, and the
two inertia bounds; and ``at most $k$ positive (resp.\ negative) eigenvalues'' is a
semialgebraic condition, because for a real symmetric matrix the number of positive
eigenvalues equals the number of sign changes in the coefficient sequence of its
characteristic polynomial (Descartes' rule is exact for polynomials with only real
roots), which is a condition on the signs of finitely many polynomials in the entries.

\emph{Step 4: closures.} Let $\pi$ be stable and put
\[
F:=\{\pl\in\Xss:\ p_i=p_j\text{ whenever }i\sim_\pi j\}
=\bigsqcup_{\pi'\succeq\pi}\big(X_{\pi'}\cap\Xss\big),
\]
a closed $\Mob$-invariant subset of $\Xss$. We claim that
$\overline{q(X_\pi)}=q(F)$.

First, $q(F)\subseteq\overline{q(X_\pi)}$: every $\ql\in F$ is a limit of
configurations in $X_\pi$ (split each block of $\pi(\ql)$ into the blocks of $\pi$ it
contains and move these to distinct points near the common value), and $q$ is
continuous.

Second, $\overline{q(X_\pi)}\subseteq q(F)$. Write $\Gamma$ also for the Gauss map
$\Polys\to X$ on polygons, $x\mapsto([w_1],\dots,[w_\ell])$; it is continuous,
$\Isom(\HH^{n+1})$-equivariant, and takes values in $\Xss$ by
Theorem~\ref{T:gauss}(2). Put $P_\pi:=\Gamma^{-1}(q(X_\pi))\subset\Po_\ell(s)$.
Since $\Xst$ is $q$-saturated and $q$ separates stable orbits (proof of
Corollary~\ref{C:compact}), $[x]\in P_\pi$ iff $q(\Gamma(x))=q(\pl)$ for some $\pl\in X_\pi$, iff
$\Gamma(x)\in\Mob X_\pi=X_\pi$; so $P_\pi$ is the image in $\Po_\ell(s)$ of
$\{x\in\Polys:\Gamma(x)\in X_\pi\}$, which is contained in the closed invariant set
$\{x\in\Polys:\Gamma(x)\in F\}$. The image of a closed invariant set under the orbit
map $\Polys\to\Po_\ell(s)$ is closed, since the preimage of its image is the set
itself. Hence $\overline{P_\pi}\subseteq\{[x]:\Gamma(x)\in F\}$, and applying the
homeomorphism $\Gamma\colon\Po_\ell(s)\to\MopsL$, under which $[x]\mapsto q(\Gamma(x))$,
gives $\overline{q(X_\pi)}=\Gamma(\overline{P_\pi})\subseteq q(F)$.

Finally, $q(F)$ is as stated in (2c). A configuration in $F$ has pattern
$\pi'\succeq\pi$. If $\pi'$ is stable it contributes $X_{\pi'}/\Mob$. Otherwise $\pi'$
has a block $A$ with $|A|=\ell/2$ (and none larger, by semistability), $A$ is a union
of blocks of $\pi$, and the configuration is sent to $\pt_{A|[\ell]\setminus A}$ by
Lemma~\ref{L:orbits}(3). Conversely, every splitting $\{A,B\}$ with $A$ a union of blocks of $\pi$
arises in this way from $\pi'=\{A,B\}$. This proves (2c) and completes the proof.
\end{proof}

The polystable quotient is thus stratified by the stable partitions, together with
finitely many points for $\ell$ even: it records which punctures coincided and where
the distinct values are, and nothing else.
\subsection{Masses from the polygon}

\begin{Proposition}\label{P:massespoly}
Assume that the Gauss points $p_1,\dots,p_\ell$ are pairwise
distinct. Then $A_{ij}:=-2G_{ij}>0$ for $i\ne j$. Let $u$ solve
$\sum_{j\ne i}(u_i+u_j+\ln A_{ij})=0$,
i.e.\ $u_i=\frac{1}{\ell-2}\big(\frac{R}{\ell-1}-R_i\big)$ with
$R_i=\sum_{j\ne i}\ln A_{ij}$, $R=\tfrac12\sum_iR_i$. Then
\begin{equation}\label{eq:masses_poly}
M_{ij}=e^{u_i+u_j}A_{ij},\qquad \eta_{ij}=-\tfrac{n-2}{2}\ln M_{ij},\qquad
\hm_i=\sum_{j\ne i}M_{ij}^{-\frac{n-2}{2}} .
\end{equation}
Here $M$ is the unique diagonal rescaling of the matrix $A =(A_{ij})$
with all row products equal to $1$.  Moreover, with
$D_{ab}=d(x_a,x_b)$ the vertex distances,
\begin{equation}\label{eq:Aij}
A_{ij}=2\big[\cosh D_{i+1,j+1}-\lambda(\cosh D_{i,j+1}+\cosh D_{i+1,j})
+\lambda^2\cosh D_{ij}\big],
\end{equation}
so, with the side length fixed, all the masses are explicit functions of the pairwise vertex distances of the polygon.
\end{Proposition}

\begin{proof}
We have $A_{ij}=t_it_jd_{ij}^2$ by \eqref{eq:chordal}, while
$d_{ij}^2=e^{2(f_i+f_j)/\kp}M_{ij}$ by the proof of Lemma~\ref{L:gram}. Taking
logarithms entrywise off the diagonal (with zero diagonal), $\ln A$
and $\ln M=-2\eta/\kp$ differ by $\varpi(\ln t+2f/\kp)\in\Lambda_\ell$.
Since $\ln M\in\Hh=\Lambda_\ell^{\perp}$, Lemma~\ref{L:proj} gives
$\ln M=P(\ln A)=\ln A+\varpi(u)$ with $u$ the unique vector making
the row sums vanish, and the computation of Lemma~\ref{L:centering}
(with $-\kp L$ replaced by $\ln A$) gives the displayed $u$. For
\eqref{eq:Aij} expand $-2\langle
x_{i+1}-\lambda x_i,x_{j+1}-\lambda x_j\rangle$ using $\langle
x_a,x_b\rangle=-\cosh D_{ab}$; for $j=i$ the bracket is
$1-2\lambda\cosh s+\lambda^2=0$, as it must be.
\end{proof}

The standard bending systems \cite{KM96,KMT} preserve the diagonal
lengths in a chosen compatible collection (equivalently, the action
variables attached to a triangulation), not all pairwise vertex
distances occurring in \eqref{eq:Aij}. Thus \eqref{eq:Aij} does not
imply that the ADM masses are preserved by bending. 

\begin{Example}[the minimizer is the regular polygon]\label{E:regular}
For $\ell\le n+2$ the regular simplex with $\ell$ vertices has all $d_{ij}$ equal,
so by symmetry the scalings solving \eqref{eq:closing} are equal and
the polygon is regular: equilateral with all diagonals of a given type
equal.
\end{Example}

\begin{Remark}\label{R:s}
Because $\Gamma:\Po_\ell(s)\xrightarrow{\ \cong\ }
\overline{\Mo}^{\,\mathrm{ps}}_\ell(S^n)$ is a homeomorphism onto the
same polystable (equivalently, closed-orbit semistable) quotient for
every $s$, the homeomorphism type of $\Po_\ell(s)$ does not depend on
$s$. In the model of Theorem~\ref{T:model} this says that moving the
affine slice inside the fixed inertia cone does not change the
topology of the intersection.
\end{Remark}

By Corollary~\ref{C:simplex} and Example~\ref{E:regular}, for $\ell\le n+2$ the
masses are minimized, simultaneously at all ends, at the regular polygon, where
$\m_i=\frac{4(\ell-1)}{n-2}$ for all $i$.

\section{Small \texorpdfstring{$\ell$}{l} and the shape-space models}\label{s:small}
\subsection{\texorpdfstring{$\ell=4$}{l=4}}
Put $P_1=d_{12}d_{34}$, $P_2=d_{13}d_{24}$, $P_3=d_{14}d_{23}$,
$u_a=\ln P_a$. The $P_a$ themselves are not conformal invariants, but
their ratios are, and $\eta$ depends only on the ratios.
\begin{Proposition}\label{P:l4}
For $\ell=4$,
\[
\eta_{12}=\eta_{34}=\tfrac{\kp}{6}(u_2+u_3-2u_1),\quad
\eta_{13}=\eta_{24}=\tfrac{\kp}{6}(u_1+u_3-2u_2),
\]
\[
\eta_{14}=\eta_{23}=\tfrac{\kp}{6}(u_1+u_2-2u_3),
\]
so with $Q=(P_1P_2P_3)^{1/3}$ one has $\xi_{ij}=(Q/P_a)^{\kp/2}$ for
the pairing $a$ containing $\{i,j\}$ and
\begin{equation}\label{eq:l4_mass}
\hm_1=\hm_2=\hm_3=\hm_4
=\Big(\frac{Q}{P_1}\Big)^{\frac{n-2}{2}}+\Big(\frac{Q}{P_2}\Big)^{\frac{n-2}{2}}
+\Big(\frac{Q}{P_3}\Big)^{\frac{n-2}{2}} .
\end{equation}
In the conformally invariant coordinates $x=P_2/P_1$, $y=P_3/P_1$,
\[
\Mo_4(S^n)\ \cong\ \{(x,y)\in\R^2_{>0}:\ |x-y|\le1\le x+y\},
\]
and $\Xi$ is the restriction to this region of the map $(x,y)\mapsto\eta$
given by the formulas above with $u=(0,\ln x,\ln y)$. The three
boundary lines $P_a=P_b+P_c$ consist of the concyclic configurations,
one line for each cyclic order. The region has three ends, one for
each split $2{+}2$; hence $\Mo_4(S^n)$ is homeomorphic to the closed
disc $D^2$ with three punctures on the boundary circle.
\end{Proposition}

\begin{proof}
The formulas follow from \eqref{eq:f_sol} by direct computation; for
instance $f_1+f_2=\frac\kp2(T_1+T_2-\frac23T)$ gives
$\eta_{12}=\frac\kp6\big[3(L_{13}+L_{14}+L_{23}+L_{24})-2T\big]
=\frac\kp6(u_2+u_3-2u_1)$, and the rest follows by symmetry. Any four
points of $S^n$ lie on a round $S^2\subset S^n$; normalizing
$p_1=0$, $p_2=1$, $p_3=\infty$ on $S^2=\C\cup\{\infty\}$ and
$p_4=z$, one finds $(x,y)=(|1-z|,|z|)$. Ptolemy's inequality gives
$|x-y|\le1\le x+y$, with equality in one of the three inequalities iff
$z$ is real, i.e.\ iff the four points are concyclic, the three cases
being the three cyclic orders. Conversely the circles $|1-z|=x$ and
$|z|=y$ meet iff $|x-y|\le1\le x+y$, so every point of the region
is attained; and $z,\bar z$ give the same $(x,y)$ and are exchanged by
the reflection in the real circle, which lies in $\Mob(S^n)$. Since
$\Xi$ is injective (Theorem~\ref{T:image}) and depends only on
$(x,y)$, the map $[\pl]\mapsto(x,y)$ is a bijection of $\Mo_4(S^n)$
onto the region, and it is a homeomorphism because $\Xi$ is one and
$(x,y)\mapsto\eta$ is a homeomorphism of $\R^2_{>0}$ onto $\mathcal
H_4\cong\R^2$. The two corners $(0,1)$, $(1,0)$ of the region and its
end along $|x-y|\le1$, $x,y\to\infty$ correspond to the collisions
$P_2\to0$, $P_3\to0$ and $P_1\to0$ respectively, i.e.\ to the three
splits $2{+}2$; adding these three points to the region gives a closed
disc.
\end{proof}
\subsection{\texorpdfstring{$\ell=5$}{l=5}}
Five points span at most a round $S^3$, so $\Mo_5(S^n)$ is the same
space for all $n\ge3$ and has dimension $5=\dim\mathcal H_5$.

\begin{Theorem}\label{T:l5}
$\Xi$ is a homeomorphism of $\Mo_5(S^n)$ onto
\[
\mathcal R_5=\{\eta\in\mathcal H_5:\ M(\eta)\ \text{has exactly one
  positive eigenvalue}\},
\]
a closed $5$-dimensional region. Its interior, $\det M>0$, consists of
the configurations spanning $S^3$; its boundary, $\det M=0$, of those
lying on a round $S^2$. The bound $\hm_i\ge4$ is attained exactly at
the regular $4$-simplex, all ten distances being $\sqrt{5/2}$.
\end{Theorem}

\begin{proof}
Theorem~\ref{T:image} with $\ell=5$: the bound on negative eigenvalues
is vacuous.  By Lemma~\ref{L:gram}, $\rank M$ is the dimension of the
span of the null vectors $v_i$, which is $r+2$ for $S^r$ the smallest
round subsphere containing the punctures (the $v_i$ then span a copy
of $\R^{1,r+1}$); so $\det M>0$ exactly for configurations spanning
$S^3$ (signature $(1,4)$) and $\det M=0$ exactly for those on a round
$S^2$ (signature $(1,3)$ and one zero eigenvalue). Signature $(1,4)$ is an open condition,
so $\{\det M>0\}\cap\mathcal R_5$ is open in $\Hh$ and consists of
interior points.

For the boundary we work in $\Mo_5(S^3)$. This is legitimate: every
five points of $S^n$ lie on a round $S^3$, so the inclusion
$S^3\subset S^n$ induces a bijection $\Mo_5(S^3)\to\Mo_5(S^n)$
compatible with the mass maps, and both mass maps are homeomorphisms
onto $\mathcal R_5$ by Theorem~\ref{T:image}; hence the bijection is
a homeomorphism. Let first $[\pl]$ be a configuration on a round
$S^2\subset S^3$ which is not concyclic. Its stabilizer $\Gamma$ in
$\Mob(S^3)$ is $\{1,\rho\}$, where $\rho$ is the reflection in the
$S^2$: an element of $\Mob(S^2)$ fixing three points is the identity
or the reflection in the circle through them, and the latter moves
every point off that circle. The action of $\Mob(S^3)$ on
$\Conf_5(S^3)$ is proper, so by the slice theorem a neighbourhood of
$[\pl]$ in $\Mo_5(S^3)$ is $U/\Gamma$ with $U\subset\R^5$ a
$\Gamma$-invariant neighbourhood of $0$ in the normal space to the
orbit. The reflection $\rho$ fixes pointwise the $4$-dimensional locus
of configurations on that $S^2$ and acts by $-1$ on the normal
direction, so $\det(\rho|_{\R^5})=-1$. By the transfer,
$H_5(U/\Gamma,U/\Gamma\setminus0;\Q)=H_5(U,U\setminus0;\Q)^\Gamma=0$.
Since $\Xi$ is a homeomorphism onto $\mathcal R_5$, the local homology
of $\mathcal R_5$ at $\Xi[\pl]$ vanishes, whereas an interior point of
a closed subset of $\R^5$ has local homology $\Q$. Hence $\Xi[\pl]$
lies in the boundary. Finally, a concyclic configuration, whose
stabilizer contains an $\mathrm{O}(2)$ and is not finite, is a limit
of non-concyclic configurations on the same $S^2$, and
$\partial\mathcal R_5$ is closed; so these lie in the boundary as well.
The last statement is Corollary~\ref{C:simplex}.
\end{proof}

The three ends of $\Mo_4$, the degeneration rays, and the fan they
generate for $\ell=5$ --- the Petersen graph --- are studied in \cite{PartII}.

\subsection{Shape-space models and contractibility}
Sending $p_\ell\mapsto\infty$ identifies the stabilizer with $\Sim(\R^n)$, so
\begin{equation}\label{eq:shape1}
\Mo_\ell(S^n)\cong\Conf_{\ell-1}(\R^n)/\Sim(\R^n)
\cong\big(S^{\,n(\ell-2)-1}\setminus\Delta\big)/O(n),
\end{equation}
Kendall's shape space of $\ell-1$ points in $\R^n$ (taken modulo
$\mathrm{O}(n)$ rather than $\mathrm{SO}(n)$, i.e.\ the reflection
shape space) with its collision locus removed.  In Gram coordinates
$Q_{ij}=\langle z_i,z_j\rangle$, $1\le i,j\le\ell-1$, for the centred
points,
\begin{equation}\label{eq:shape2}
\Mo_\ell(S^n)\cong\{Q\succeq0,\ Q\mathbf1=0,\ \operatorname{tr}Q=1,\ \rank Q\le n\}
\setminus\bigcup_{i<j}\{Q_{ii}+Q_{jj}-2Q_{ij}=0\},
\end{equation}
where $\mathbf1\in\R^{\ell-1}$ is the all-ones vector: a rank-truncated
spectrahedron minus some of its exposed faces; for $\ell\le n+2$ the
rank condition is vacuous, since $\rank Q\le\ell-2$. The third model
is the image $\mathcal R^{(n)}_\ell:=\Xi(\Mo_\ell(S^n))\subset\Hh$, described by
Theorem~\ref{T:image}(2), in
which $\rank M(\eta)=r+2$ for $S^r$ the spanning subsphere.

\begin{Proposition}\label{P:contract}
For $3\le\ell\le n+2$, $\Mo_\ell(S^n)$ is contractible, and its image
$\mathcal R^{(n)}_\ell$ is star-shaped about $\eta=0$.
\end{Proposition}

\begin{proof}
For $\ell\le n+2$ the rank condition in \eqref{eq:shape2} is vacuous,
so the moduli space $\Mo_\ell(S^n)\cong K\setminus F$ where
$K=\{Q\succeq0,\ Q\mathbf 1=0,\ \operatorname{tr}Q=1\}$ is
compact convex and $F$ a union of exposed faces,
$Q_{ii}+Q_{jj}-2Q_{ij}=(e_i-e_j)^TQ(e_i-e_j)\ge0$ on $K$. The Gram
matrix $Q_0$ of the regular simplex lies in the relative interior of
$K$ and not in $F$, and for every $Q\in K$ the half-open segment
$[Q_0,Q)$ lies in the relative interior of $K$, hence misses $F$; so
  $K\setminus F$ is star-shaped about $Q_0$, and contractible. This is
  also in \cite{BS1}.

For the second statement: the space of matrices $\Hh$ is linear, so
$t\eta\in\Hh$ for $t\in[0,1]$. For $0<t\le1$, $M(t\eta)=M(\eta)^{\circ
  t}=\Lambda^tD^{\circ t}\Lambda^t$, where $D=(d_{ij}^2)$ and
$M(\eta)=\Lambda D\Lambda$. Now $D$ is a squared Euclidean distance
matrix, equivalently $-D$ is conditionally positive semidefinite on
the hyperplane $\mathbf 1^\perp$. By Schoenberg's theorem \cite{Sch},
$D^{\circ t}=(d_{ij}^{2t})$ is again a squared Euclidean distance
matrix. Hence it has exactly one positive eigenvalue and at most
$\ell-1$ negative eigenvalues. Since $\ell-1\le n+1$, the conditions
of Theorem~\ref{T:image}(2) hold for $M(t\eta)$, inertia being
preserved by the positive diagonal congruence. At $t=0$, $M(0)=J-I$,
which has inertia $(1,\ell-1)$ and therefore also satisfies the same
condition.
\end{proof}

\begin{Remark}\label{R:snow}
For $\ell>n+2$ the origin is not in the image
(Corollary~\ref{C:simplex}), so the radial argument has no base point.
Moreover, snowflaking $D\mapsto D^{\circ t}$ may raise the Euclidean
embedding dimension, so Schoenberg's theorem no longer keeps
$M(t\eta)$ within the bound of $n+1$ negative eigenvalues. Whether the
image is star-shaped about some other point is open. The first case
not covered is $\ell=n+3$, which is the critical case resolved in
\cite{BS1}.
\end{Remark}
\begin{Remark}
Topology of the moduli spaces $\mathcal{M}_\ell(S^n)$ and of their
compactifications is addressed in \cite{BS1,BS2}.
\end{Remark}

\end{document}